\documentclass[reqno,11pt]{amsart}

\newtheorem{theorem}{Theorem}[section]
\newtheorem{lemma}[theorem]{Lemma}
\newtheorem{corollary}[theorem]{Corollary}

\theoremstyle{definition}

\newtheorem{assumption}{Assumption}[section]

 \theoremstyle{remark}
\newtheorem{remark}[theorem]{Remark}

\newcommand\bH{\mathbb{H}}
\newcommand\bL{\mathbb{L}}

\newcommand\bR{\mathbb{R}}

\newcommand\cF{\mathcal{F}}

\newcommand\cR{\mathcal{R}}

\newcommand{\WO}[2]{\overset{\scriptscriptstyle0}%
{W}\,\!^{#1}_{#2}}

\newcommand{\osc}{\operatornamewithlimits{osc\,}}
\newcommand{\nliminf}{\operatornamewithlimits{\underline{lim}\,}}

\newcommand{\mysection}[1]{\section{#1}
\setcounter{equation}{0}}

\newcommand\cbrk{\text{$]$\kern-.15em$]$}} 
\newcommand\opar{
\text{\,\raise.2ex\hbox{${\scriptstyle |}$}\kern-.34em$($}} 
\newcommand\cpar{%
\text{$)$\kern-.34em\raise.2ex\hbox{${\scriptstyle |}$}}\,}
\newcommand\obrk{\text{$[$\kern-.15em$[$}}

\begin{document}

\title[Maximum principle]
{Maximum principle for SPDEs 
and its applications}
\author{N.V. Krylov}
\address{127 Vincent Hall, University of Minnesota, Minneapolis,
 MN, 55455}
\thanks{The work   was partially supported by
NSF Grant DMS-0140405}
\email{krylov@math.umn.edu}
 \keywords{Maximum principle, H\"older continuity,
stochastic partial differential equations}

\renewcommand{\subjclassname}{%
\textup{2000} Mathematics Subject Classification}

\subjclass{60H15, 35R60}

\begin{abstract}
The maximum principle for SPDEs is established
in   multidimensional $C^{1}$ domains. An application
is given to proving the H\"older
continuity  up to the boundary of solutions of one-dimensional
SPDEs.
\end{abstract}

\maketitle

The maximum principle is one of the most powerful
tools in the theory of second-order
elliptic and parabolic
partial differential equations. However, until now
it did not play any significant role in the theory
of SPDEs. In this paper we show how to
apply it to   one-dimensional SPDEs
on the half line $\bR_{+}=(0,\infty)$ and prove the H\"older
continuity of solutions on $[0,\infty)$.
This result was previously known when
the coefficients of the first order derivatives of solution
appearing in the stochastic term in the equation
obeys a quite unpleasant condition. On the other hand,
if they just vanish, then the H\"older continuity
was well known before (see, for instance,
\cite{KyK1} and the references therein).

To the best of our knowledge the maximum principle was first
proved in \cite{KR81} (see also \cite{KR82}
for the case of random coefficients) for SPDEs
in the whole space by the method of random characteristics
introduced there and also in   \cite{Ku}.
Later the method of random characteristics was used
in many papers for various purposes, for instance,
to prove smoothness of solutions
(see, for instance, \cite{BG1}, \cite{BG2},
\cite{PT}, \cite{Tu} and the references therein). 
It was very tempting
to try to use this method for proving the maximum
principle for SPDEs in domains.
However, the implementation of the method
 turns out to become extremely
cumbersome and
inconvenient if the coefficients of the equation
are random processes. Also, it  requires
more regularity of solutions than actually needed.

Here in Section \ref{section 1.14.1}
 we state the maximum principle in domains
under minimal assumptions. We prove it
in Section \ref{section 2.28.2} by using methods taken
from PDEs after we prepare some auxiliary results
in Section \ref{section 2.28.1}.
 
Section \ref{section 2.22.2} contains
an application of the maximum principle
 to investigating
the H\"older continuity up to the boundary
of solutions of one-dimensional SPDEs.
Note that, for instance,
 in \cite{BG1}, \cite{BG2} and in many other
papers that can be found
from our list of references
the regularity properties are proved only inside domains.
Quite sharp regularity for solutions of SPDEs in
multidimensional domains
is established in \cite{KyK}, it is stated in terms
of appropriate weighted 
Sobolev spaces and, unfortunately, do not imply
even the pointwise continuity up to the boundary. It is worth saying
that we only deal with one-dimensional case and coefficients
independent of the space variable. In a subsequent
paper we intend to treat the general case.
In Section \ref{section 2.22.1} we introduce
some auxiliary functions used in 
Section \ref{section 2.22.2}.

We denote by $\bR^{d}$ the Euclidean space of points
$x=(x^{1},...,x^{d})$,
$$
D_{i}=\frac{\partial}{\partial x^{i}}.
$$
For a domain $D\subset\bR^{d}$ and  we set $W^{1}_{2}(D)$
to be the closure of the set of infinitely differentiable
functions $\phi$ having finite norm
$$
\|\phi\|_{W^{1}_{2}(D)}^{2}=\|\phi\|^{2}_{L_{2}(D)}
+\|\phi_{x}\|^{2}_{L_{2}(D)}
$$
with respect to this norm. Here $\phi_{x}$ is the gradient of
$\phi$. By $\WO{1}{2}(D)$ we denote the closure
of $C_{0}^{\infty}(D)$ with respect to the norm
$\|\cdot\|_{W^{1}_{2}(D)}$. 
Our way to say that $u\leq v$ on $\partial D$ is
that $(u-v)^{+}\in\WO{1}{2}(D)$.
As usual, the summation convention is
enforced and writing $N(....)$ is to say that the constant $N$ depends
and depends only on the contents 
of the parentheses. Such constants may change from line to line.

Few typos in the original version of the article
were kindly pointed out by Kyeong-Hun Kim.
The author is sincerely grateful for that.
 
\mysection{The maximum principle}

                                         \label{section 1.14.1}

Let $D$ be a domain in $\bR^{d}$ of class $C^{1}_{loc}$ and let
 $(\Omega,\cF,P)$ be
a complete probability space with
a given filtration  $(\cF_{t},t\geq0)$ of
$\sigma$-fields $\cF_{t}\subset\cF$ complete
with respect to $\cF,P$.

We are investigating some properties of a function
$u_{t}(x)=u_{t}(\omega,x)$ satisfying
$$
(\phi,u_{t })=(\phi,u_{0})+
\int_{0}^{t}
(\phi,\sigma^{ik}_{s}D_{i}u_{s}+\nu^{k}_{s}u_{s}+g^{k}_{s})
 \,dm^{k}_{s}
$$
\begin{equation}
                                               \label{1.4.1}
+\int_{0}^{t}(\phi,D_{i}(a^{ij}_{s}D_{j}u_{s})
+b^{i}_{s}D_{i}u_{s}+D_{i}(a ^{i}_{s}u_{s})-c_{s}u_{s} +f_{s}
+D_{i}f^{i}_{s})
\,dV_{s}.
\end{equation}
for all $t\in[0,\infty)$ and any $\phi\in C^{\infty}_{0}(D)$.
Here $m^{k}_{t}$, $k=1,2,...$, are  
one-dimensional continuous local $\cF_{t}$-martingales, 
starting at zero, $V_{t}$ is a nondecreasing continuous
$\cF_{t}$-adapted process starting at zero,
  $(\phi,\cdot)$ is the pairing between
a generalized function on $D$ and a test function $\phi$,
the summation convention over repeated
indices is enforced,
and the meaning of the remaining objects
and further assumptions  are
described below.  We need some real-valued
functions $\xi^{i}_{t}(x)$,
$K_{1}(t)>0$, and $K_{2}(t)\geq0$
defined for $i=1,...,d$,  $t\in[0,\infty)$, $x\in\bR^{d}$
and also depending on $\omega$. 

We assume that $a^{ij}_{t}(x)$, $b^{i}_{t}(x)$, $a ^{i}_{t}(x)$,
$c_{t}(x)$,
$\sigma^{ik}_{t}(x)$, $\nu^{k}_{t}(x)$, and $g^{k}_{s}$ 
 are real-valued  
functions defined for $i,j=1,...,d$, $k=1,2,...$,
$t\in[0,\infty)$, $x\in\bR^{d}$ and also depending on
$\omega\in\Omega$.
 
\begin{assumption}
                                          \label{assumption 2.15.1}
We suppose that, for any $\omega$,
  $ \langle m^{i},
m^{j}\rangle_{t}=0$ if $i\ne j$, and for any $k$
we have $d\langle m^{k}\rangle_{t}\leq dV_{t}$.
\end{assumption}

\begin{assumption}
                                          \label{assumption 1.10.1} 
For all values of the arguments

(i)   $\sigma^{i}:=(\sigma^{i1},\sigma^{i2},...)$,
$\nu:=(\nu^{1},\nu^{2},...)$,
$g:=(g^{1},g^{2},...)\in\ell_{2}$;

(ii) for all $\lambda\in\bR^{d}$ 
$$
|\sum_{i}\lambda^{i}\xi^{i}|^{2}\leq
K_{1}(2a^{ij}-\alpha^{ij})\lambda^{i}\lambda^{j},
$$
 where
$\alpha^{ij}=(\sigma^{i},\sigma^{j})_{\ell_{2}}$.

The case $\xi\equiv0$ is not excluded and in this case
Assumption \ref{assumption 1.10.1} (ii) is just the
usual parabolicity assumption.

\end{assumption}
\begin{assumption}
                                          \label{assumption 1.12.1}

(i) The functions
$a^{ij}_{t}(x)$, $b^{i}_{t}(x)$, $a ^{i}_{t}(x)$,
$c_{t}(x)$,
$\sigma^{ik}_{t}(x)$, $\nu^{k}_{t}(x)$,
 $\xi^{i}_{t}(x)$, 
$K_{1}(t)$, and $K_{2}(t)$ are measurable
with respect to
$(\omega,t,x)$ and $\cF_{t}$-adapted for each~$x$;

(ii)  the functions
$a^{ij}_{t}(x)$, $b^{i}_{t}(x)$, $a ^{i}_{t}(x)$,
$c_{t}(x)$,
$\sigma^{ik}_{t}(x)$, $\nu^{k}_{t}(x)$, and $\xi^{i}_{t}(x)$
 are bounded;

(iii) for each $\omega,t$
the functions
$$
\eta^{i}_{t}:=a ^{i}_{t}-b^{i}_{t}
 -(\sigma^{i}_{t},\nu_{t})_{\ell_{2}}
-\xi^{i}_{t}
$$
are once continuously differentiable on $D$,
have bounded derivatives, and satisfy
\begin{equation}
                                                \label{1.10.2}
D_{i}\eta^{i}-2c+|\nu|^{2}_{\ell_{2}}\leq K_{2}
\end{equation}
for all values of arguments;

(iv) for each $\phi\in C^{\infty}_{0}(D)$
the processes $\phi f_{t}$, $\phi f^{1}_{t}$,...,
$\phi f^{d}_{t}$  are   $L_{2}(D)$-valued and $\phi g_{t}$
is an $L_{2}(D,\ell_{2})$-valued
$\cF_{t}$-adapted and jointly measurable;
 for all $t\in[0,\infty)$ and $\omega\in\Omega$
$$
\int_{0}^{t}\big(\|\phi f_{s}\|^{2}_{L_{2}(D)}+\sum_{i}
\|\phi f^{i}_{s}\|^{2}_{L_{2}(D)}+
\|\phi g_{s}\|^{2}_{L_{2}(D,\ell_{2})}
+K_{1}(s)+K_{2}(s)\big)\,dV_{s}<\infty.
$$

\end{assumption}

\begin{assumption}
                                \label{assumption 1.15.1} 
For each $\phi\in C^{\infty}_{0}(D)$

(i) the process
$\phi u_{t}=\phi u_{t}(\omega )$ is 
  $L_{2}(D)$-valued, $\cF_{t}$-adapted,
and jointly measurable;

(ii) for any
 $\omega$
$$
\phi u_{t}\in W^{1}_{2}(D)\quad\text{($dV_{t}$-a.e.)};
$$

(iii) for each $t\in[0,\infty)$ and $\omega$
$$
\int_{0}^{t}\|\phi u_{s}\|^{2}_{W^{1}_{2}(D)}\,dV_{s}<\infty.
$$

\end{assumption}

 The above assumptions are supposed to hold throughout this section.
Here is the maximum principle saying, in particular, that if
$g^{k}=f^{i}=0$, $f\leq0$
and $u\leq0$ on the parabolic boundary of $[0,T]\times D$,
then $u\leq0$ in $[0,T] $. By the way, our solutions are 
$L_{2,loc}(D)$-valued functions
of $\omega$ and $t$, so that for
each $\omega$ and $t$ an equivalence class
is specified. 
Naturally, if we write $u_{t}(\omega)\leq0$, or $u_{t}\leq0$
we mean that in the corresponding class there is a
nonpositive function. 
\begin{theorem}
                                          \label{theorem 1.14.1}

Let $\tau_{2}\geq\tau_{1}$ be  stopping times,
$\tau_{1}<\infty$ for any $\omega$. Suppose that,
 for any $\omega$, $i=1,...,d$, $k=1,2,...$,
$$
I_{u_{t}>0}g^{k}_{t}=I_{u_{t}>0}f^{i}_{t}=0,\quad
 u^{+}_{t}
\in\WO{1}{2}(D),\quad I_{u_{t}>0}f_{t}\leq0
$$
$dV_{t}$-almost everywhere on
$(\tau_{1},\tau_{2}) $ and suppose that
$u_{\tau_{1}}\leq0$  for any $\omega$. 
 Then almost surely
$u_{t}\leq0$ for all $t\in[\tau_{1},
\tau_{2}]\cap[\tau_{1},\infty)$.
\end{theorem}
 
The following comparison principle is a generalization
of Theorem \ref{theorem 1.14.1}.

\begin{theorem}
                                          \label{theorem 3.18.1}

Let $\tau_{2}\geq\tau_{1}$ be  stopping times,
$\tau_{1}<\infty$ for any $\omega$.
 Let $\rho_{t}\geq0$, $t\in[0,\infty)$, be a  
nondecreasing continuous
$\cF_{t}$-adapted process and let
$\bar{f}_{t}$, $\bar{f}^{1}_{t}$,...,
$\bar{f}^{d}_{t}$, and $\bar{g}_{t}$ satisfy
Assumption \ref{assumption 1.12.1} (iv). Let
$\bar{u}_{t}$ be a process satisfying
Assumption \ref{assumption 1.15.1} and such that
equation \eqref{1.4.1} holds for all $t\in[0,\infty)$ 
and any $\phi\in C^{\infty}_{0}(D)$ with
$\bar{f}_{t}$, $\bar{f}^{1}_{t}$,...,
$\bar{f}^{d}_{t}$, and $\bar{g}_{t}$
in place of  $f_{t}$, $f^{1}_{t}$,...,
$f^{d}_{t}$, and $g_{t}$, respectively.

Assume that, for any $\omega$,   
($dV_{t}$-a.e.) on $[\tau_{1},\tau_{2}] $ we have
$$
I_{u_{t}>\rho\bar{u}_{t}}(g_{t}-\rho_{t}\bar{g}_{t})=I_{u_{t}>
\rho_{t}\bar{u}_{t}} (f^{i}_{t}-\rho_{t}\bar{f}^{i}_{t} )
=0, \quad i=1,...,d,
$$ 
 $$
I_{u_{t}>\rho_{t}\bar{u}_{t}}(f _{t}-\rho_{t}
\bar{f}_{t})\leq 0,\quad
I_{u_{t}>\rho_{t}\bar{u}_{t}}\bar{u}_{t}\geq0,\quad
 (u_{t}-\rho_{t}\bar{u}_{t})^{+} 
\in\WO{1}{2}(D).
$$
Finally, assume that
$u_{\tau_{1}}\leq\rho_{\tau_{1}}\bar{u}_{\tau_{1}}$
   for any $\omega$. 
 
Then almost surely
$u_{t}\leq\rho_{t}\bar{u}_{t}$ for all $t\in[\tau_{1},
\tau_{2}]\cap[\tau_{1},\infty)$.
\end{theorem}

\begin{corollary}
                                        \label{corollary 2.10.1}
Assume that, for any $\omega$,   
($dV_{t}$-a.e.) on $(\tau_{1},\tau_{2})\times D$ we have
$$
I_{u_{t}>1}(\nu^{k}_{t}+g^{k}_{t})=I_{u_{t}>1}(f^{i}_{t}
+a ^{i}_{t})=0, \quad i=1,...,d,k=1,2,...
$$ 
 $$
I_{u_{t}>1}f\leq I_{u_{t}>1}c,\quad  (u_{t}-1)^{+} 
\in\WO{1}{2}(D).
$$
Also assume that
$u_{\tau_{1}}\leq1$   for any $\omega$. Then almost surely
$u_{t}\leq1$ for all $t\in[\tau_{1},
\tau_{2}]\cap[\tau_{1},\infty)$.
\end{corollary}

Indeed,  it suffices to take $\bar{u}_{t}\equiv1$,
$\rho_{t}\equiv1$ and observe that $\bar{u}_{t}$
satisfies \eqref{1.4.1} with $\bar{f}^{i}_{t}=-a^{i}$, $\bar{f}=c$,
and $\bar{g}=-\nu_{t} $ in place of $f^{i}_{t}$, $f_{t}$,
and $g_{t}$,
respectively.

This corollary generalizes the corresponding results
of \cite{KR81} and \cite{KR82}, where $\nu^{k}=g^{k}=
f^{i}=a^{i}=0$.

\begin{remark}
                                     \label{remark 2.20.1}
Our equation has a special structure, which may look
quite restrictive.
In particular, we assume that the martingales
$m^{k}_{t}$ are mutually orthogonal. The general case,
actually, reduces to this particular one after 
using the fact that
one can always orthogonalize the martingales
by using, for instance, the Gramm-Schmidt procedure.
This, of course, would change $\sigma$, $\nu$, and $g$,
and writing the corresponding general conditions
would only obscure the matter. Then passing from
$m^{k}_{t}$ to (no summation in $k$)
$$
\int_{0}^{t}\rho^{k}_{s}\,dm^{k}_{s},\quad
\rho^{k}_{s}=\big(\frac{dt}{dt+d\langle m^{k}\rangle_{t}}
\big)^{1/2}
$$
allows one to have $d\langle m^{k}\rangle_{t}\leq dt$
and adding after that  $t$ to $V_{t}$ allows one
to have $d\langle m^{k}\rangle_{t}\leq dV_{t}$.
Again we should modify our coefficients but
we will see in the proof of 
Theorem \ref{theorem 3.18.1} that this modification
does not affect Assumption \ref{assumption 1.10.1},
  which is an assumption about parabolicity
of our equation and not strict nondegeneracy.

\end{remark}

\mysection{Auxiliary results}
                                   \label{section 2.28.1}

In this section the notation $u_{t}$
is sometimes used for different objects than
in Section \ref{section 1.14.1}.

Denote by $\cR$ the set of real-valued functions convex
$r(x)$ on $\bR$ such that

\noindent(i) $r$ is continuously differentiable, $r(0)=r'(0)=0$,

\noindent(ii) $r'$ is absolutely continuous, its
derivative $r''$ is bounded and left continuous,
that is usual $r''$ which exists almost everywhere
is bounded and there is a left-continuous
function with which $r''$ coincides
almost everywhere.

For $r\in\cR$ by $r''$ we will always mean the left-continuous
modification of the usual second-order derivative of $r$.

\begin{remark}
                                         \label{remark 1.8.1}
For each $r\in\cR$ there exists a sequence
$r_{n}\in\cR$ of {\em infinitely\/}
 differentiable functions
 such that $|r_{n}(x)|\leq N|x|^{2}$,
$|r'_{n}(x)|\leq N|x|$, and $|r''_{n}|\leq N$ with $N
<\infty$ independent of
$x$ and $n$,
$r_{n},r_{n}',r_{n}''\to r,r'r''$ on $\bR$.

Indeed, let $\zeta\in C^{\infty}_{0}(\bR)$ be a nonnegative
function with support in $(0,1)$ and unit integral. For 
$\varepsilon>0$ define $\zeta_{\varepsilon}(x)
=\varepsilon^{-1}\zeta(x/\varepsilon)$ and
$r_{\varepsilon }(x)=r*\zeta_{\varepsilon}(x)
-r*\zeta_{\varepsilon}(0)-xr'*\zeta_{\varepsilon}(0)$.
Then $r_{\varepsilon}$ is infinitely differentiable,
$r_{\varepsilon}(0)=r'_{\varepsilon}(0)=0$,
$$
|r''_{\varepsilon}|=
|r''*\zeta_{\varepsilon}|\leq\sup|r''|<\infty.
$$
In particular,
$$
|r'_{\varepsilon}(x)|=|\int_{0}^{x}r''_{\varepsilon}
(y)\,dy|\leq N|x|,\quad
|r_{\varepsilon}(x)|=|\int_{0}^{x}r'_{\varepsilon}
(y)\,dy|\leq N|x|^{2}.
$$
Finally, the convergences $r_{\varepsilon}\to r$
and $r'_{\varepsilon}\to r'$ follow by
the continuity
of $r$ and $r'$ and the convergence $r''_{\varepsilon}
\to r''$ follows from the dominated convergence theorem,
the left continuity of $r''$ and the formula
$$
r''_{\varepsilon}(x)=\int_{0}^{1}r''(x-\varepsilon y)
\zeta(y)\,dy.
$$

\end{remark}

In the following lemma the assumption that $D$
is a locally smooth domain is
not used.
\begin{lemma}
                                        \label{lemma 12.19.1}

Let  
$u_{t}=u_{t}(\omega )$ be  an 
  $L_{2}(D)$-valued process
 such that $u_{0}$ is $\cF_{0}$-measurable.
Let $f_{t}$ and $g_{t}=(g^{1}_{t},g^{2}_{t},...)$
be $\cF_{t}$-adapted and jointly measurable
processes with values
$L_{2}(D)$ and $L_{2}(D,\ell_{2})$, respectively.
Assume that  for each $t\in[0,\infty)$  we have
\begin{equation}
                                               \label{12.20.1}
\int_{0}^{t}(\|f_{s}\|^{2}_{L_{2}(D)}+
\|g_{s}\|^{2}_{L_{2}(D,\ell_{2})})\,dV_{s}<\infty
\end{equation}
and for any $\phi\in C^{\infty}_{0}(D)$
\begin{equation}
                                               \label{12.19.1}
(\phi,u_{t})_{L_{2}(D)}=(\phi,u_{0})_{L_{2}(D)}+
\int_{0}^{t}(\phi,f_{s})_{L_{2}(D)}\,dV_{s}+\int_{0}^{t}
(\phi,g^{k}_{s})_{L_{2}(D)}\,dm^{k}_{s}.
\end{equation}

Then  (i) $u_{t}$ is a continuous $L_{2}(D)$-valued
function (a.s.); 
(ii) for any $r\in\cR$ (a.s.)
for all $t\in[0,\infty)$ 
\begin{equation}
                                               \label{12.19.02}
\|r^{1/2}(u_{t})\|^{2}_{L_{2}(D)}=
\|r^{1/2}(u_{0})\|^{2}_{L_{2}(D)}+
\int_{0}^{t}h_{s}\,dV_{s}+m_{t},
\end{equation}
where
\begin{equation}
                                               \label{1.8.2}
h_{s}:=(r'(u _{s}),f_{s})_{L_{2}(D)}
+(1/2)\|(r'')^{1/2}(u _{s}) \check{g}_{s}\|^{2}_{L_{2}(D,\ell_{2})},
\end{equation}
$$
\check{g}_{s}^{k}:=\bigg(\frac{d\langle m^{k}\rangle_{s}}
{dV_{s}}\bigg)^{1/2}g^{k}_{s},\quad
m_{t}:=\int_{0}^{t}
(r'(u_{s}),g^{k}_{s})_{L_{2}(D)}\,dm^{k}_{s}
$$
and $m_{t}$ is a local martinagale;

(iii) (a.s.)
for $t\in[0,\infty)$ 
\begin{equation}
                                               \label{12.19.2}
\|u^{+}_{t}\|^{2}_{L_{2}(D)}=
\|u^{+}_{0}\|^{2}_{L_{2}(D)}+
\int_{0}^{t}h_{s}\,dV_{s}+m_{t},
\end{equation}
where
$$
h_{s}:=2(u^{+}_{s},f_{s})_{L_{2}(D)}
+\|\check{g}_{s}I_{u_{s}>0}\|^{2}_{L_{2}(D,\ell_{2})},
$$
$$
m_{t}:=2\int_{0}^{t}
(u^{+}_{s},g^{k}_{s})_{L_{2}(D)}\,dm^{k}_{s}
$$
and $m_{t}$ is a local martinagale.

\end{lemma}

Proof. (i) Recall that the operation of stochastic integration
of Hilbert space valued processes is well defined.
Therefore, the process
$$
\hat{u}_{t}=u_{0}+\int_{0}^{t}f_{s}\,dV_{s}+\int_{0}^{t}
g^{k}_{s}\,dm^{k}_{s}
$$
is well defined as a continuous $L_{2}(D)$-valued
process. We also recall how the scalar product interacts with
integrals. Then it is seen that for any $t$ and $\phi\in
C^{\infty}_{0}(D)$ we have $(\phi,u_{t})=(\phi,\hat{u}_{t})$
(a.s.). Since both parts are continuous in $t$, the equality holds
for all $t$ at once (a.s.), and since $C^{\infty}_{0}(D)$
is dense in $L_{2}(D)$, we have that $u_{t}=\hat{u}_{t}$
for all $t$ (a.s.). This proves (i). As a corollary
we obtain that
\begin{equation}
                                               \label{12.19.3}
\sup_{t\leq T}\|u_{t}\|_{L_{2}(D)}<\infty,
\quad\forall T<\infty\quad\text{(a.s.)}.
\end{equation}

(ii) It suffices to prove \eqref{12.19.02} for
infinitely differentiable $r\in\cR$. Indeed,
for $r_{n}$
from Remark \ref{remark 1.8.1}, passing to the limit
in all term in \eqref{12.19.02} apart from $m_{t}$
presents no problem at all in light of \eqref{12.19.3}
and the dominated convergence theorem.
Also
$$
m _{t}(n):=\int_{0}^{t}
(r_{n}'(u_{s}),g^{k}_{s})_{L_{2}(D)}\,dm^{k}_{s}\to m_{t}
$$
uniformly in $t$ on finite intervals in probability
because
$$
\langle m(n)-m \rangle_{t}=\int_{0}^{t}
\sum_{k}\big((r'_{n}-r')(u_{s}),g^{k}_{s})_{L_{2}(D)}^{2}
\,d\langle m^{k}\rangle_{s}
$$
$$
\leq \int_{0}^{t}
\sum_{k}\big((r'_{n}-r')(u_{s}),g^{k}_{s})_{L_{2}(D)}^{2}
\,dV_{s}
$$
$$
\leq \int_{0}^{t}\|(r'_{n}-r')(u_{s})\|^{2}
_{L_{2}(D)} \|
g_{s}\|^{2}_{L_{2}(D,\ell_{2})}\,dV_{s} \to0
$$
again owing to \eqref{12.19.3}
and the dominated convergence theorem.

Thus, we may concentrate on the case that 
$r$ is infinitely differentiable.
Take a symmetric $\zeta\in C^{\infty}_{0}(\bR^{d})$
with support in the unit ball centered at the origin
and unit integral. For $\varepsilon>0$ set $\zeta_{\varepsilon}(x)
=\varepsilon^{-d}\zeta(x/\varepsilon)$ and for functions $v=v(x)$
define $v^{(\varepsilon)}=v*\zeta_{\varepsilon}$. Also set
$$
D_{\varepsilon}=\{x\in D:\text{dist}\,(x,\partial D)<\varepsilon
\}.
$$
 According to \eqref{12.19.1}  
for any $x\in D_{\varepsilon}$ and $t\geq0$
we have
$$
u^{(\varepsilon)}_{t}(x)=u^{(\varepsilon)}_{0}(x)
+\int_{0}^{t}f^{(\varepsilon)}_{s}(x)\,dV_{s}
+\int_{0}^{t}g^{k(\varepsilon)}_{s}(x)\,dm^{k}_{s}.
$$

By It\^o's formula
we have  that   on $ D_{\varepsilon}$
$$
r(u^{(\varepsilon)}_{t}) =
r(u^{(\varepsilon)}_{0})
+\int_{0}^{t}[r '(u^{(\varepsilon)}_{s})f^{(\varepsilon)}_{s}
+(1/2)r '' (u^{(\varepsilon)}_{s})
|\check{g}^{(\varepsilon)}_{s}|^{2}_{\ell_{2}}]\,dV_{s}
$$
\begin{equation}
                                               \label{12.20.2}
+\int_{0}^{t}r ' (u^{(\varepsilon)}_{s})
g^{k(\varepsilon)}_{s}\,dm^{k}_{s}.
\end{equation}
Here, for each $\varepsilon>0$,
 the integrands are smooth functions of $x$ 
and their magnitudes along with the 
magnitudes of each of their
derivatives in $D_{\varepsilon}$ are
majorated by
a constant (possibly depending on $\varepsilon$)
times
$$
\|f_{s}\|_{L_{2}(D)}\sup_{s\leq t}\|u_{s}\|_{L_{2}(D)}
\quad \text{or}\quad \|g_{s}\|^{2}_{L_{2}(D,\ell_{2})}
\quad \text{or}\quad 
\|g^{k}_{s}\|_{L_{2}(D)}\sup_{s\leq t}\|u_{s}\|_{L_{2}(D)}.
$$
 This and \eqref{12.20.1} and \eqref{12.19.3} allow us to use Fubini's
theorem while integrating through \eqref{12.20.2}
 and conclude
$$
\|r ^{1/2} (u^{(\varepsilon)}_{t})\|^{2}
_{L_{2}(D_{\varepsilon})}
=\|r^{1/2} (u^{(\varepsilon)}_{0})\|^{2}
_{L_{2}(D_{\varepsilon})}+
\int_{0}^{t}[(r ' (u^{(\varepsilon)}_{s}),
f^{(\varepsilon)}_{s})_{L_{2}(D_{\varepsilon})}
$$
$$
+(1/2)\|(r '' )^{1/2}(u^{(\varepsilon)}_{s})
 \check{g}^{(\varepsilon)}_{s}\|^{2}_{L_{2}(D_{\varepsilon},\ell_{2})}]
\,dV_{s}
$$
\begin{equation}
                                               \label{12.20.3}
+\int_{0}^{t}(r ' (u^{(\varepsilon)}_{s}),
g^{k(\varepsilon)}_{s})_{L_{2}(D_{\varepsilon})}\,dm^{k}_{s}.
\end{equation}

Now we let $\varepsilon\downarrow0$. We use that
for any function $v\in L_{2}(D)$
$$
\|v^{(\varepsilon)}\|_{L_{2}(D_{\varepsilon})}
\leq
\|(vI_{D})^{(\varepsilon)}\|_{L_{2}(\bR^{d})}
\leq\|vI_{D}\|_{L_{2}(\bR^{d})}=\|v\|_{L_{2}(D)}
$$
and $v^{(\varepsilon)}I_{D_{\varepsilon}}\to v$
in $L_{2}(D)$. In particular,
$g^{k(\varepsilon)}_{s}I_{D_{\varepsilon}}\to
g^{k }_{s}$ and
$  u^{(\varepsilon)}_{s}I_{D_{\varepsilon}} \to
 u _{s} $ implying that 
$r '(u^{(\varepsilon)}_{s})I_{D_{\varepsilon}}\to
r ' (u _{s})$ in $L_{2}(D)$
 and 
$$
(r ' (u^{(\varepsilon)}_{s}),
g^{k(\varepsilon)}_{s})_{L_{2}(D_{\varepsilon})}
\to (r ' (u _{s}),
g^{k }_{s})_{L_{2}(D )}
$$
for each $k$ and $dP\times dV_{s}$-almost all $(\omega,s)$.

We also use \eqref{12.20.1} and \eqref{12.19.3}
to assert that
$$
\sum_{k=1}^{\infty}
\int_{0}^{t}\sup_{\varepsilon\in(0,1)}
|(r '(u^{(\varepsilon)}_{s}),
g^{k(\varepsilon)}_{s})_{L_{2}(D_{\varepsilon})}|^{2}\,dV_{s}
$$
$$
\leq N\sup_{s\leq t}\|u_{s}\|^{2}_{L_{2}(D)}
\int_{0}^{t}\|
g_{s}\|_{L_{2}(D,\ell_{2})}^{2}\,dV_{s}<\infty.
$$
As is easy to see this implies that the local martingale part in
\eqref{12.20.3} converges  to $m_{t}$ as $\varepsilon\downarrow0$
in probability locally uniformly with respect
to $t$.

Similar manipulations with other terms in
\eqref{12.20.3} allow us to get
 \eqref{12.19.02}. Since \eqref{12.19.2}
is just a particular case of \eqref{12.19.02}, the lemma is proved.

\begin{remark}
Lemma \ref{lemma 12.19.1} remains true if in the definition
of $\cR$
 instead
of requiring $r''$ to have a left-continuous
modification we required it to have a right-continuous
one, and of course, in \eqref{1.8.2} used this
right-continuous
modification. This is seen after replacing $u$ with $-u$.

In case $r(x)=(x^{+})^{2}$ the function $r''$
has both right- and left-continuous modifications,
so that in the definition of $m_{t}$ one can use
$2I_{u_{s}>0}$ or $2I_{u_{s}\geq0}$. It follows that (a.s.)
for any $t$
$$
\int_{0}^{t}\| \check{g}_{s}I_{u_{s}=0}\|^{2}_{L_{2}(D,\ell_{2})}\,dV_{s}
=0.
$$
Furthermore, since, for any $v\in L_{2}(D)$, $u_{t}+v$ has the same
form as $u_{t}$,
$$
\int_{0}^{t}\|\check{g}_{s}I_{u_{s}=v}\|^{2}_{L_{2}(D,\ell_{2})}\,dV_{s}
=0.
$$
\end{remark}

\begin{lemma}
                                        \label{lemma 3.18.2}
Let $D$ be an arbitrary  domain. Let
$ u_{t} $ be  an 
  $L_{2}(D)$-valued 
$\cF_{0}$-measurable process such that for any~$\omega$
$$
 u_{t}\in W^{1}_{2}(D)
$$
($dV_{t}$-a.e.) and for each $T\in[0,\infty)$ and $\omega$
\begin{equation}
                                               \label{3.18.4}
\int_{0}^{T}\| u_{t}\|^{2}_{W^{1}_{2}(D)}\,dV_{t}<\infty.
\end{equation}
Let $f_{t},f^{1}_{t},...,f^{d}_{t}$, and
$g_{t}=(g^{1}_{t},g^{2}_{t},...)$ be $\cF_{t}$-adapted and jointly
measurable processes with values
in $L_{2}(D)$ and $L_{2}(D,\ell_{2})$, respectively.
Assume that  for each $t\in[0,\infty)$  we have
\begin{equation}
                                               \label{3.18.5} 
\int_{0}^{t}(\| f_{s}\|^{2}_{L_{2}(D)}+\sum_{i}
\| f^{i}_{s}\|^{2}_{L_{2}(D)}+
\| g_{s}\|^{2}_{L_{2}(D,\ell_{2})})\,dV_{s}<\infty,
\end{equation}
and  for each $t\in[0,\infty)$, $\phi\in C^{\infty}_{0}(D)$,
and $\omega$
\begin{equation}
                                               \label{3.18.6}
(\phi,u_{t})=(\phi,u_{0})+
\int_{0}^{t}(\phi,f_{s}+D_{i}f^{i}_{s})\,dV_{s}+\int_{0}^{t}
(\phi,g^{k}_{s})\,dm^{k}_{s}.
\end{equation}
Finally, assume that 
there is a compact set $G\subset D$
such that 
$$
u_{t}(x)=f_{t}(x)=f^{i}_{t}(x)
=g^{k}_{t}(x)=0
$$
 outside $G$. Then 

 (a) $u_{t}$ is a continuous $L_{2}(D)$-valued
function (a.s.); (b) (a.s.) for all $t\in[0,\infty)$
\begin{equation}
                                                  \label{1.8.8}
\|u^{+}_{t}\|^{2}_{L_{2}(D)}=
\|u^{+}_{0}\|^{2}_{L_{2}(D)}+
\int_{0}^{t}h _{s}\,dV_{s}+m _{t} ,
\end{equation}
where
$$
h _{s}:=2( u^{+}_{s}, f_{s})_{L_{2}(D)}
-2(I_{u_{s}>0}D_{i} u_{s}, f^{i}_{s})_{L_{2}(D)}
+\|
\check{g}_{s}I_{\phi u_{s}>0}\|^{2}_{L_{2}(D,\ell_{2})},
$$
$$
m_{t} :=2\int_{0}^{t}
( u^{+}_{s}, g^{k}_{s})_{L_{2}(D)}\,dm^{k}_{s}.
$$
\end{lemma}

Proof. Observe that \eqref{3.18.6} holds for all
infinitely differentiable functions $\phi$. Furthermore,
since $u_{t}\in \WO{1}{2}(D)$ ($dV_t$-a.e.)
assertion (a) is well known  (see, for instance,
\cite{KR}, the references therein, and Remark \ref{remark 2.20.1}).

To prove (b),
take $\varepsilon$ smaller than the distance between $G$
and $\partial D$. Notice that, owing to the symmetry
of $\zeta$, for $\phi\in C^{\infty}_{0}(D)$
$$
(\phi^{(\varepsilon)},u_{t})=(\phi^{(\varepsilon)},u_{t})
_{L_{2}(\bR^{d})}=
(\phi,u^{(\varepsilon)}_{t})
_{L_{2}(\bR^{d})}=
(\phi,u^{(\varepsilon)}_{t})
_{L_{2}(D)}.
$$
Therefore, it follows from \eqref{3.18.6} that
for any $\phi\in C^{\infty}_{0}(D)$
\begin{equation}
                                               \label{1.5.7}
(\phi, u_{t} ^{(\varepsilon)})
=(\phi, u_{0} ^{(\varepsilon)})
+\int_{0}^{t}(\phi, \tilde{f}^{\varepsilon}_{s})\,ds
+\int_{0}^{t}(\phi,  g^{k(\varepsilon)}_{s} )\,dm^{k}_{s},
\end{equation}
where
$$
\tilde{f}^{\varepsilon}_{s}:= 
f^{(\varepsilon)}_{s} +D_{i} 
f^{i(\varepsilon)}_{s}  
$$
is an $L_{2}(D)$-valued function with norm that is
locally square integrable
against $dV_{s}$. By  Lemma \ref{lemma 12.19.1} for any
$r\in\cR$
\begin{equation}
                                               \label{1.5.8}
\|r^{1/2}( u_{t} ^{(\varepsilon)})\|^{2}_{L_{2}(D )}
=\|r^{1/2}( u_{0} ^{(\varepsilon)})\|^{2}_{L_{2}(D )}
+\int_{0}^{t}h^{\varepsilon}_{s}\,dV_{s}+m^{\varepsilon}_{t},
\end{equation}
where 
$$
m^{\varepsilon}_{t}:=\int_{0}^{t}
(r'( u_{s} ^{(\varepsilon)}),
\check{g}^{k(\varepsilon)}_{s} )_{L_{2}(D )}\,dm^{k}_{s},
$$
$$
h^{\varepsilon}_{s}:=
(r'( u_{s} ^{(\varepsilon)}),\tilde{f}^{\varepsilon}_{s})
_{L_{2}(D)}+(1/2)\|(r'')^{1/2}( u_{s} ^{(\varepsilon)})  \check{g}^{
(\varepsilon)}_{s}  
\|^{2}_{L_{2}(D ,\ell_{2})}
$$
$$
=(r'( u_{s} ^{(\varepsilon)}),f^{(\varepsilon)}_{s})_{L_{2}(D)}
-(r''( u_{s} ^{(\varepsilon)})
D_{i} u_{s} ^{(\varepsilon)},f^{i(\varepsilon)}_{s})
_{L_{2}(D)}+(1/2)\|(r'')^{1/2}( u_{s} ^{(\varepsilon)})  \check{g}^{
(\varepsilon)}_{s}  
\|^{2}_{L_{2}(D ,\ell_{2})}.
$$

If $r$ is infinitely differentiable, then
by using \eqref{3.18.4} and \eqref{3.18.5} one easily
passes  to the limit in \eqref{1.5.8} as $
\varepsilon\to0$. The argument is quite similar
to the corresponding argument in the proof
of Lemma \ref{lemma 12.19.1} and, for
smooth $r\in\cR$, yields
\begin{equation}
                                               \label{1.5.08}
\|r^{1/2}( u_{t})\|^{2}_{L_{2}(D )}
=\|r^{1/2}( u_{0})\|^{2}_{L_{2}(D )}
+\int_{0}^{t}h_{s}\,dV_{s}+m_{t},
\end{equation}
where 
$$
m_{t}=\int_{0}^{t}
(r'( u_{s}),
g^{k}_{s} )_{L_{2}(D )}\,dm^{k}_{s},
$$
$$
h_{s}
=(r'( u_{s}),f_{s})_{L_{2}(D)}
-(r''( u_{s} )
D_{i} u_{s},f^{i}_{s})
_{L_{2}(D)}+(1/2)\|(r'')^{1/2}( u_{s})\check{g}_{s}  
\|^{2}_{L_{2}(D ,\ell_{2})}.
$$

Finally, as in the proof of  Lemma \ref{lemma 12.19.1}
one easily passes from smooth $r\in\cR$
to arbitrary ones and gets \eqref{12.19.2}
by taking $r(x)=(x^{+})^{2}$. The lemma is proved.

Lemma \ref{lemma 3.18.2} serves as an auxiliary tool to
prove a deeper result.

\begin{lemma}
                                        \label{lemma 3.18.1}
Let $D$ be an arbitrary  domain.
Assume that for each $\phi\in C^{\infty}_{0}(D)$ 

(i) 
$\phi u_{t} $ is  an 
  $L_{2}(D)$-valued process
 such that $\phi u_{0}$ is $\cF_{0}$-measurable;
  
(ii) for any~$\omega$
$$
\phi u_{t}\in W^{1}_{2}(D)
$$
($dV_{t}$-a.e.) and for each $T\in[0,\infty)$ and $\omega$
\begin{equation}
                                               \label{1.8.04}
\int_{0}^{T}\|\phi u_{t}\|^{2}_{W^{1}_{2}(D)}\,dV_{t}<\infty.
\end{equation}
(iii) Let $f_{t},f^{1}_{t},...,f^{d}_{t}$, and
$g_{t}=(g^{1}_{t},g^{2}_{t},...)$ be $\cF_{t}$-adapted and jointly
measurable processes with values
in $L_{2}(D)$ and $L_{2}(D,\ell_{2})$, respectively.
Assume that  for each $t\in[0,\infty)$ 
and $\phi\in C^{\infty}_{0}(D)$ we have
\begin{equation}
                                               \label{1.5.03} 
\int_{0}^{t}(\|\phi f_{s}\|^{2}_{L_{2}(D)}+\sum_{i}
\|\phi f^{i}_{s}\|^{2}_{L_{2}(D)}+
\|\phi g_{s}\|^{2}_{L_{2}(D,\ell_{2})})\,dV_{s}<\infty,
\end{equation}
\begin{equation}
                                               \label{1.5.04}
(\phi,u_{t})=(\phi,u_{0})+
\int_{0}^{t}(\phi,f_{s}+D_{i}f^{i}_{s})\,dV_{s}+\int_{0}^{t}
(\phi,g^{k}_{s})\,dm^{k}_{s}.
\end{equation}
 
Then, for any $\phi\in C^{\infty}_{0}(D)$,

 (a) $\phi u_{t} $ is a continuous $L_{2}(D)$-valued
function (a.s.); (b) (a.s.) for all $t\in[0,\infty)$
\begin{equation}
                                                  \label{1.8.8}
\|(\phi u_{t})^{+}\|^{2}_{L_{2}(D)}=
\|(\phi u_{0})^{+}\|^{2}_{L_{2}(D)}+
\int_{0}^{t}h _{s}\,dV_{s}+m _{t} ,
\end{equation}
where
$$
h _{s}:=2((\phi u_{s})^{+},\phi f_{s}
-f^{i}_{s}D_{i}\phi)_{L_{2}(D)}
-2(I_{\phi u_{s}>0}D_{i}(\phi u_{s}),\phi f^{i}_{s})_{L_{2}(D)}
$$
$$
+\|
\phi\check{g}_{s}I_{\phi u_{s}>0}\|^{2}_{L_{2}(D,\ell_{2})},
\quad
m_{t} :=2\int_{0}^{t}
(\phi u^{+}_{s},\phi g^{k}_{s})_{L_{2}(D)}\,dm^{k}_{s}.
$$
 \end{lemma}

Proof. Clearly, for any $\phi,\eta\in C^{\infty}_{0}(D)$
we have
\begin{equation}
                                                   \label{3.18.7}
(\phi,\eta u_{t})=(\phi,\eta u_{0})+
\int_{0}^{t}(\phi,\eta f_{s}-f^{i}_{s}D_{i}\eta
+D_{i}(\eta f^{i}_{s})
)\,dV_{s}+\int_{0}^{t}
(\phi,\eta  g^{k}_{s})\,dm^{k}_{s}.
\end{equation}

Therefore, $\eta u_{t}$ satisfies the assumptions
of Lemma \ref{lemma 3.18.2} with
$\eta f_{s}-f^{i}_{s}D_{i}\eta$,
$\eta f^{i}_{s}$, and $\eta g^{k}_{s}$
in place of $f_{s}$, $D_{i}f^{i}_{s}$, and $g^{k}_{s}$,
respectively.

By applying  Lemma \ref{lemma 3.18.2} to $\eta u_{t}$
in  place of $u_{t}$ we get the result with $\eta$
in place of $\phi$. This certainly proves the lemma.

\mysection{Proof of Theorems \protect\ref{theorem 1.14.1}
and \protect\ref{theorem 3.18.1}}
                                         \label{section 2.28.2}

In this section the assumptions stated in Section
\ref{section 1.14.1} are supposed to be 
satisfied.
We use the fact that
due to our hypothesis that $D\in C^{1}$,
there exist  sequences $\zeta_{n}$ and $\bar{\zeta}_{n}$ of
nonnegative
$C^{\infty}_{0}(D)$-functions such that
$0\leq\zeta_{n},\bar{\zeta}_{n}\leq1$, $\zeta_{n},\bar{\zeta}_{n}
\to1$ in $D$
as $n\to\infty$  and for any $v\in\WO{1}{2}(D)$, $i=1,...,d$,
\begin{equation}
                                               \label{1.8.6}
\| vD_{i}\zeta_{n}\|_{L_{2}(D)}\leq  
N(\|(1-\bar{\zeta}_{n}) v\|_{L_{2}(D)}+\|(1-
\bar{\zeta}_{n})D v\|_{L_{2}(D)}),
\end{equation}
where $N$ is independent of $n$ and $v$
(see, for
instance, the proof of Theorem 5.5.2 in \cite{Ev}). 
We also know
(see, for instance, the proof of Lemma 2.3.2 in \cite{LU}
or  Problem 17,
Chapter 5 of \cite{Ev}) that if $v\in W^{1}_{2}(D)$,
then $v^{+}\in W^{1}_{2}(D)$ and 
$$
D_{i}v^{+}=I_{v>0}
D_{i}v.
$$

{\bf Proof of Theorem \ref{theorem 1.14.1}}.
Set
$$
K=K_{1}+K_{2},\quad \varphi_{t}=\int_{0}^{t}K(s)\,ds.
$$
Take the sequences of nonnegative $\zeta_{n},\bar{\zeta}_{n}
\in C^{\infty}_{0}(D)$ from above.
By It\^o's formula and Lemma \ref{lemma 3.18.1}
$$
\|\zeta_{n} u^{+}_{t}\|^{2}_{L_{2}(D)}e^{- \varphi_{t}}
=\|\zeta_{n} u_{0}^{+}\|_{L_{2}(D)}^{2}+
\int_{0}^{t}h^{n}_{s}\,dV_{s}+m_{t}(n),
$$
where
$$
e^{ \varphi_{s}}
h^{n}_{s}=I_{1s}+I_{2s}+I_{3s}- K(s)
\|\zeta_{n} u^{+}_{s}\|^{2}_{L_{2}(D)},
$$
$$
I_{1s}=2\big(\zeta_{n}u^{+}_{s},\zeta_{n}[f_{s}+b^{i}_{s}D_{i}u_{s}
-c_{s}u_{s}]-[u_{s}a^{i}_{s}+a^{ij}_{s}D_{j}u_{s}
+f^{i}_{s}]D_{i}\zeta_{n}\big)_{L_{2}(D)} ,
$$
$$
I_{2s}=
-2\big(I_{\zeta_{n}u_{s}>0} D_{i}(\zeta_{n}u_{s}),\zeta_{n}
[u_{s}a^{i}_{s}+a^{ij}_{s}D_{j}u_{s}
+f^{i}_{s}]\big)_{L_{2}(D)},
$$
$$
I_{3s}=\|\zeta_{n}I_{\zeta_{n}u_{s}>0}
[\check{\sigma}^{i}_{s}D_{i}u_{s}+
\check{\nu}_{s}u_{s}+\check{g}_{s}]\|^{2}_{L_{2}(D,\ell_{2})},
$$
$$
m_{t}(n)=2\int_{0}^{t}e^{- \varphi_{s}}\big(\zeta_{n}u^{+}_{s},
\zeta_{n}[\sigma^{ik}_{s}D_{i}u_{s}+\nu^{k}_{s}
u_{s}+g^{k}_{s}]\big)_{L_{2}(D)}
\,dm^{k}_{s}.
$$
Since $u^{+}_{\tau_{1}}=0$ we have 
\begin{equation}
                                                 \label{2.18.1}
e^{-\varphi_{\tau_{2}\wedge t\vee\tau_{1}}}
\|\zeta_{n}u^{+}_{\tau_{2}\wedge t\vee\tau_{1}}\|_{L_{2}(D)}^{2 }
=\int_{0}^{t}I_{\tau_{2}>s>\tau_{1}}h^{n}_{s}\,dV_{s}
+\bar{m}_{t}(n),
\end{equation}
where
$$
\bar{m}_{t}(n):=
 m_{\tau_{2}\wedge t\vee\tau_{1}}(n)-m_{\tau_{1}}(n) 
$$
is a local martingale.

Next we use the assumptions of the theorem and see that
for $dV_{s}$-almost all $s\in(\tau_{1},\tau_{2})$ we have
$$
I_{1s}\leq 2(\zeta_{n}u^{+}_{s},\zeta_{n}[ b^{i}_{s}D_{i}u_{s}
-c_{s}u_{s}]-[u_{s}a^{i}_{s}+a^{ij}_{s}D_{j}u_{s}
 ]D_{i}\zeta_{n})_{L_{2}(D)} 
$$
$$
=2(\zeta_{n}^{2}u^{+}_{s},  b^{i}_{s}D_{i}u^{+}_{s}
-c_{s}u^{+}_{s})_{L_{2}(D)}+I_{4s}
$$
with
$$
I_{4s}=
-
2(\zeta_{n}u^{+}_{s}D_{i}\zeta_{n}
, u^{+}_{s}a^{i}_{s}+a^{ij}_{s}D_{j}u^{+}_{s}
  )_{L_{2}(D)} . 
$$
At this moment we recall \eqref{1.8.6}
and observe that
$$
(\zeta_{n}|D_{i}\zeta_{n}|
, (u^{+}_{s})^{2})_{L_{2}(D)}
\leq N\|\zeta_{n}u^{+}_{s}\|_{L_{2}(D)}
\|u^{+}_{s}D_{i}\zeta_{n}\|_{L_{2}(D)}.
$$
Then we see that

$$
I_{4s}\leq N\big(\|(1-\bar{\zeta}_{n}) u^{+}_{s}\|_{L_{2}(D)}+\|
(1-\bar{\zeta}_{n})D u^{+}_{s}\|_{L_{2}(D)}\big)
\|u^{+}_{s}\| _{W^{1}_{2}(D)},
$$
where 
and below by $N$ we  denote various  finite constants.

In $I_{2s}$
$$
I_{\zeta_{n}u_{s}>0} D_{i}(\zeta_{n}u_{s})=
 D_{i}(\zeta_{n}u^{+}_{s})=u^{+}_{s}D_{i}\zeta_{n}
+\zeta_{n}D_{i}u^{+}_{s},
$$
so that
$$
I_{2s}=-2(\zeta_{n}^{2}D_{i}u^{+}_{s},u^{+}_{s}a^{i}_{s}+
a^{ij}_{s}D_{j}u^{+}_{s})_{L_{2}(D)}+I_{4s}.
$$
Next,
$$
\zeta_{n}I_{\zeta_{n}u_{s}>0}=\zeta_{n}I_{ u_{s}>0},
$$
$$
I_{3s}\leq \|\zeta_{n}I_{ u_{s}>0}
[ \sigma ^{i}_{s}D_{i}u_{s}+
 \nu_{s} u_{s} ]\|^{2}_{L_{2}(D,\ell_{2})}=
(\zeta_{n}^{2}D_{i}u^{+}_{s},\alpha^{ij}_{s}D_{j}
u^{+}_{s})_{L_{2}(D)}
$$
$$
+2(\zeta_{n}^{2}D_{i}u^{+}_{s},u^{+}_{s}(\sigma^{i}_{s},
\nu_{s})_{\ell_{2}})_{L_{2}(D)}+\|\zeta_{n}
|\nu_{s}|_{\ell_{2}}u^{+}_{s}\|^{2}_{L_{2}(D)}.
$$
Also observe that certain parts of $I_{2s}$ and $I_{3s}$
can be combined if we use that
$$
-2(\zeta_{n}^{2}D_{i}u^{+}_{s}, 
a^{ij}_{s}D_{j}u^{+}_{s})_{L_{2}(D)}+
(\zeta_{n}^{2}D_{i}u^{+}_{s},\alpha^{ij}_{s}D_{j}
u^{+}_{s})_{L_{2}(D)}
$$
$$
\leq -K_{1}^{-1}(s)
\|\zeta_{n}\xi^{i}_{s}D_{i}u^{+}_{s}\|^{2}_{L_{2}(D)}.
$$
It follows that for $dV_{s}$-almost all $s\in(\tau_{1},\tau_{2})$
$$
e^{ \varphi_{s}}h^{n}_{s}\leq \int_{D}[\zeta_{n}^{2}
(b^{i}_{s}-a^{i}_{s}+(\sigma^{i}_{s},\nu_{s})_{\ell_{2}}
)2u^{+}_{s}D_{i}u^{+}_{s}+\zeta_{n}^{2}
(u^{+}_{s})^{2}(|\nu_{s}|^{2}_{\ell_{2}}-2c_{s})]\,dx
$$
$$
+N\big(\|(1-\bar{\zeta}_{n}) u^{+}_{s}\|_{L_{2}(D)}+\|(1-
\bar{\zeta}_{n})D u^{+}_{s}\|_{L_{2}(D)}\big)
\|u^{+}_{s}\| _{W^{1}_{2}(D)}
$$
$$
-K_{1}^{-1}(s)
\|\zeta_{n}\xi^{i}_{s}D_{i}u^{+}_{s}\|^{2}_{L_{2}(D)}-K(s)
\|\zeta_{n} u^{+}_{s}\|^{2}_{L_{2}(D)}.
$$
Here
$$
b^{i}_{s}-a^{i}_{s}+(\sigma^{i}_{s},\nu_{s})_{\ell_{2}}
=-\xi^{i}_{s}-\eta^{i}_{s}
$$
and we transform the integral of
$$
\zeta_{n}^{2}(-\eta^{i}_{s})2u^{+}_{s}D_{i}u^{+}_{s}
=-\eta^{i}_{s}\zeta^{2}_{n}D_{i}(u^{+}_{s})^{2}
$$
by integrating by parts. Then we get
that for $dV_{s}$-almost all $s\in(\tau_{1},\tau_{2})$
$$
e^{\varphi_{s}}h^{n}_{s}\leq \int_{D}[-\zeta_{n}^{2}
\xi^{i}_{s}2u^{+}_{s}D_{i}u^{+}_{s}+
\zeta_{n}^{2}(u^{+}_{s})^{2}
(|\nu_{s}|^{2}_{\ell_{2}}-2c_{s}+D_{i}\eta^{i}_{s})]\,dx
$$
$$
+2\int_{D}(u^{+}_{s})^{2}\zeta_{n}\eta^{i}_{s}D_{i}\zeta_{n}\,dx
+N\big(\|(1-\bar{\zeta}_{n}) u^{+}_{s}\|_{L_{2}(D)}+\|(1-
\bar{\zeta}_{n})D u^{+}_{s}\|_{L_{2}(D)}\big)
\|u^{+}_{s}\| _{W^{1}_{2}(D)}
$$
$$
-K_{1}^{-1}(s)
\|\zeta_{n}\xi^{i}_{s}D_{i}u^{+}_{s}\|^{2}_{L_{2}(D)}-K(s)
\|\zeta_{n} u^{+}_{s}\|^{2}_{L_{2}(D)}.
$$

We also use the fact that
$$
|-\zeta_{n}^{2}\xi^{i}_{s} 2u^{+}_{s}D_{i}u^{+}_{s}|
\leq K_{1}^{-1}(s)|\zeta_{n}\xi^{i}_{s}D_{i}u^{+}_{s}|^{2}+
K_{1}(s)\zeta_{n}^{2}(u^{+}_{s})^{2},
$$
$$
|\nu_{s}|^{2}_{\ell_{2}}-2c_{s}+D_{i}\eta^{i}_{s}
\leq K_{2}(s).
$$
Then we easily see that
for $dV_{s}$-almost all $s\in(\tau_{1},\tau_{2})$
$$
e^{\varphi_{s}}h^{n}_{s}\leq 
 N\big(\|(1-\bar{\zeta}_{n}) u^{+}_{s}\|_{L_{2}(D)}+\|
(1-\bar{\zeta}_{n})D u^{+}_{s}\|_{L_{2}(D)}\big)
\|u^{+}_{s}\| _{W^{1}_{2}(D)}
$$
$$
+2\int_{D}(u^{+}_{s})^{2}\zeta_{n}\eta^{i}_{s}D_{i}\zeta_{n}\,dx
\leq  N\big(\|(1-\bar{\zeta}_{n}) u^{+}_{s}\|_{L_{2}(D)}+\|
(1-\bar{\zeta}_{n})D u^{+}_{s}\|_{L_{2}(D)}\big)
\|u^{+}_{s}\| _{W^{1}_{2}(D)}.
$$
Now \eqref{2.18.1} yields
$$
e^{-\varphi_{\tau_{2}\wedge t\vee\tau_{1}}}
\|u^{+}_{\tau_{2}\wedge t\vee\tau_{1}}\|_{L_{2}(D)}^{2 }
\leq \bar{m}_{t}(n)
$$
\begin{equation}
                                                 \label{2.19.1}
+N
\int_{0}^{t}(\|(1-\bar{\zeta}_{n}) u^{+}_{s}\|_{L_{2}(D)}+\|
(1-\bar{\zeta}_{n})D u^{+}_{s}\|_{L_{2}(D)})
\|u^{+}_{s}\| _{W^{1}_{2}(D)}\,dV_{s} .
\end{equation}

The integrals against $dV_{s}$ in \eqref{2.19.1}
tend to zero as $n\to\infty$ by the dominated
convergence theorem.
 Since the sum of them with   continuous local martingales
is nonnegative, the local martingales and the right-hand
side of \eqref{2.19.1} tend to zero uniformly
on finite time intervals in probability
(see, for instance, \cite{Kr95}).
So does the left-hand side and 
the theorem is proved.

{\bf Proof of Theorem \ref{theorem 3.18.1}}.
Obviously, $\hat{u}_{t}=\rho_{t}\bar{u}_{t}$ satisfies
$$
(\phi,\hat{u}_{t })=(\phi,\hat{u}_{0})+
\int_{0}^{t}
(\phi,\sigma^{ik}_{s}D_{i}\hat{u}_{s}+\nu^{k}_{s}\hat{u}_{s}+
\rho_{s}\bar{g}^{k}_{s})
 \,dm^{k}_{s}+\int_{0}^{t}(\phi,\bar{u}_{s})\,d\rho_{s}+
$$
\begin{equation}
                                                   \label{2.19.2}
+\int_{0}^{t}(\phi,D_{i}(a^{ij}_{s}D_{j}\hat{u}_{s})
+b^{i}_{s}D_{i}\hat{u}_{s}+
D_{i}(a ^{i}_{s}\hat{u}_{s})-c_{s}\hat{u}_{s} 
+\rho_{s}\bar{f}_{s}
+\rho_{s}D_{i}\bar{f}^{i}_{s})
\,dV_{s}.
\end{equation}
We rewrite this equation introducing 
$$
\hat{V}_{t}=V_{t}+\rho_{t},\quad p_{t}=\frac{d \rho_{t}}
{d\hat{V}_{t}},\quad
q_{t}=\frac{dV_{t}}
{d\hat{V}_{t}},
$$
$$
(\hat{a}^{ij}_{t},\hat{a}^{i}_{t},
\hat{b}^{i}_{t},\hat{c}_{t})=q_{t}
( a^{ij}_{t}, a^{i}_{t},
 b^{i}_{t}, c_{t}),\quad  
(\hat{\sigma}^{ik}_{t},\hat{\nu}^{k}_{t})=q^{1/2}_{t}
(\sigma^{ik}_{t},\nu^{k}_{t}),
$$
$$
\hat{f}_{t}=q_{t}\rho_{t}\bar{f}_{t}+p_{t}\bar{u}_{t},\quad
\hat{f}^{i}_{t}=q_{t}\rho_{t}\bar{f}^{i}_{t},\quad
\hat{g}^{k}_{t}=q_{t}^{1/2}\rho_{t}\bar{g}^{k}_{t}.
$$
 
We also set
$$
\hat{m}_{t}^{k}=\int_{0}^{t}q_{s}^{-1/2}\,dm^{k}_{s}
\quad\quad(0^{-1/2}:=0).
$$
Notice that since $d\langle m^{k}\rangle_{t}
\leq dV_{t}=q_{t}d\hat{V}_{t}$ the last integral
makes sense.

In this notation \eqref{1.4.1} and \eqref{2.19.2}
are rewritten as
$$
(\phi,u_{t })=(\phi,u_{0})+
\int_{0}^{t}
(\phi,\hat\sigma^{ik}_{s}D_{i}u_{s}+
\hat\nu^{k}_{s}u_{s}+q^{1/2}_{s}g^{k}_{s})
 \,d\hat m^{k}_{s}
$$
$$
+\int_{0}^{t}(\phi,D_{i}(\hat a^{ij}_{s}D_{j}u_{s})
+\hat b^{i}_{s}D_{i}u_{s}
+D_{i}(\hat a ^{i}_{s}u_{s})-\hat c_{s}u_{s} +q_{s}f_{s}
+q_{s}D_{i}f^{i}_{s})
\,d\hat V_{s},
$$
$$
(\phi,\hat u_{t })=(\phi,\hat u_{0})+
\int_{0}^{t}
(\phi,\hat\sigma^{ik}_{s}D_{i}\hat u_{s}+
\hat\nu^{k}_{s}\hat u_{s}+\hat g^{k}_{s})
 \,d\hat m^{k}_{s}
$$
$$
+\int_{0}^{t}(\phi,D_{i}(\hat a^{ij}_{s}D_{j}\hat u_{s})
+\hat b^{i}_{s}D_{i}\hat u_{s}
+D_{i}(\hat a ^{i}_{s}\hat u_{s})-\hat c_{s}\hat u_{s} +\hat f_{s}
+  D_{i}\hat f^{i}_{s})
\,d\hat V_{s},
$$
respectively. We subtract these equations, denote
$v_{t}=u_{t}-\hat{u}_{t}$, and observe that
for any $\omega$ we have $d\hat V_{s}$-almost everywhere
on $(\tau_{1},\tau_{2})$ that
$$
I_{v_{s}>0}(q^{1/2}_{s}g^{k}_{s}-\hat g^{k}_{s})
=I_{v_{s}>0}(q _{s}f^{i}_{s}-\hat{f}^{i}_{s})=0,
$$
$$
 I_{v_{s}>0}(q _{s}f _{s}-\hat{f} _{s})=
q _{s} I_{v_{s}>0}(f _{s}-\rho_{s}\bar{f} _{s})
-p_{s}\bar{u}_{s}I_{v_{s}>0}\leq0.
$$
We also use the fact that the above versions of
equations \eqref{1.4.1} and \eqref{2.19.2}
satisfy the same Assumptions \ref{assumption 2.15.1},
\ref{assumption 1.10.1}, \ref{assumption 1.12.1}, and
\ref{assumption 1.15.1} with $q _{s}\xi^{i}_{s}$ and
$q _{s}K_{i}(s)$ 
in place of $\xi^{i}_{s}$ and $K_{i}(s)$, respectively.
Then we the desired result directly from
Theorem \ref{theorem 1.14.1}. The theorem is proved.

\mysection{Auxiliary functions}
                                         \label{section 2.22.1}

Let $C[0,\infty)$ be the set of real-valued continuous functions
on $[0,\infty)$. For $x_{\cdot}\in C[0,\infty)$ 
set $x_{s}=x_{0}$ for $s\leq0$ and for  $n=0,1,2,...$ and $t\geq0$
introduce
$$
\Delta_{n}^{-}(x_{\cdot},t)=2^{n/2}
\osc_{[t-2^{-n},t ]} x_{\cdot}\,.
$$
If $c\in(0,\infty)$, then define
$$
M_{n}^{-}(x_{\cdot},c,t)=\#\{k=0,...,n: 
\Delta^{-}_{k}(x_{\cdot},t)\leq  c\}.
$$
For $n$ negative we set $M_{n}^{-}(x_{\cdot},c,t):=0$.
For $c\geq0,d>0,\delta>0$ introduce
$$
\gamma(c,d,\delta)=1-P(\min_{t\leq\delta/2}w_{t}\leq-c-
d/\sqrt{2},\max_{t\leq\delta/2}w_{t}\leq d-d/\sqrt{2}).
$$
As is easy to see
$$
\gamma(c,d,\delta)\geq P( w_{t}\quad
\text{reaches}\quad d-d/\sqrt{2}\quad
\text{before reaching}\quad -c-d/\sqrt{2})
$$
$$
=\frac{c+d/\sqrt{2}}{c+d}>1/\sqrt{2},
$$
so that $2\log_{2}\gamma(c,d,\delta)>-1 $.

Set, for $m=0, 1, 2,...$,
$$
 Q_{m}:=Q_{m}(x_{\cdot}):=\{(s,y):s\geq0, x_{s}<y<x_{s}+2^{-m/2}\}.
$$

\begin{lemma}
                                         \label{lemma 1.25.1}
For $m=0,1,2,...$, $t\geq 0$ and $x\in(0,2^{-m/2})$ introduce
$$
r_{m}(t,x)=
r_{m}(x_{\cdot},t,x)=
P(x_{t}+x+w_{\tau}\sqrt{\delta}=x_{t-\tau}+2^{-m/2}),
$$
where $\tau=\inf\{s>0:
(t-s,x_{t}+x+w_{s}\sqrt{\delta})\not\in Q_{m}\}$. Then
\begin{equation}
                                              \label{1.25.1}
r_{m}(t,x)\leq
[\gamma(c,d,\delta)]^{M_{m+n}^{-}(x_{\cdot},c,t)
-M_{m-1}^{-}(x_{\cdot},c,t)-k},
\end{equation}
where
$n=n(2^{ m/2}x/d)$, $k=k(c+d)$, and
 $$
n(y)=[(-2\log_{2}y)_{+}]
,\quad k(d)=2+[(2\log_{2} d )_{+}].
$$

\end{lemma}

Proof. Define
$$
\bar{t}=2^{m }t,\quad \bar{x}=2^{m/2}x,\quad
\bar{w}_{s}=2^{m/2}w_{s2^{-m}},\quad
 \bar{x}_{s}=2^{m/2}x_{s2^{-m}}.
$$
Then as is easy to see $r_{m}(t,x)$
is rewritten as
\begin{equation}
                                              \label{1.25.2}
P(\bar{x}_{\bar{t}}+\bar{x}
+\bar{w}_{\bar{\tau}}\sqrt{\delta}=
\bar{x}_{\bar{t}-\bar{\tau}}+1),
\end{equation}
where
$$
\bar{\tau}=\inf\{s>0:
(\bar{t}-s,\bar{x}_{\bar{t}}+\bar{x}+
\bar{w}_{s}\sqrt{\delta})\not\in Q_{0}(\bar{x}_{\cdot})\}
=2^{m}\tau.
$$
Since $\bar{w}_{\cdot}$ is a Wiener process,
by Corollary 3.4 of \cite{Kr03} expression \eqref{1.25.2}
is less than
$$
[\gamma(c,d,\delta)]^{M_{\bar{n}}^{-}(\bar{x}_{\cdot},c,\bar{t})-k},
$$
where $\bar{n}=n(\bar{x}/d)$. Here
$$
M_{\bar{n}}^{-}(\bar{x}_{\cdot},c,\bar{t})
=\#\{j=0,...,\bar{n}:2^{j/2}\osc_{[\,\bar{t}-2^{-j},\bar{t}\,]}
\bar{x}_{\cdot}\leq c\}
$$
$$
=\#\{j=0,...,\bar{n}:2^{(j+m)/2}\osc_{[t-2^{-j-m},t\,]}
x_{\cdot}\leq c\}
$$
$$
=\#\{j=m,...,m+\bar{n}:\Delta^{-}_{j}(x_{\cdot})\leq c\}
=M_{m+\bar{n}}^{-}(x_{\cdot},c,t)
-M_{m-1}^{-}(x_{\cdot},c,t) 
$$
and the result follows. The lemma is proved.

\begin{lemma}
                                         \label{lemma 1.25.2}
Let $T\in(0,\infty)$. Assume that
$$
\nliminf_{m\to\infty}\frac{1}{m+1}
\inf_{t\in[0,T]}M^{-}_{m}(x_{\cdot},c,t)>\alpha>0.
$$
Take constants $p>0$ and $\nu$ so that
\begin{equation}
                                               \label{2.7.2}
 1<\nu p<p\chi+1<0 ,
\end{equation}
where $\chi=-2\alpha\log_{2}\gamma(c,d,\delta)$.
Then, for $r_{m}$ from Lemma \ref{lemma 1.25.1} it
holds that
\begin{equation}
                                              \label{1.25.3}
\sup_{m\geq0}\sup_{t\in[0,T]}2^{-m(\nu p-1)/(2\alpha)}
\int_{0}^{2^{- m/2}}\frac{1}{x^{\nu p}}r^{p}_{m}(t,x)
\,dx<\infty.
\end{equation}

\end{lemma}

Proof. By Lemma \ref{lemma 1.25.1} for a constant $N$
and $\gamma=\gamma(c,d,\delta)$
$$
r_{m}(t,x)\leq N\gamma^{M^{-}_{m+n}-m},
$$
where $x\leq 2^{-m/2}$, $n=n(2^{m/2}x/d)$,
$M^{-}_{m+n}=M^{-}_{m+n}(x_{\cdot},c,t)$.
Furthermore,
$$
n(2^{m/2}x/d)\geq (-2\log_{2}(2^{m/2}x)+2\log_{2}d)_{+}-1
$$
$$
\leq  (-2\log_{2}(2^{m/2}x))_{+}-N=-m-2\log_{2}x-N,
$$
where $N$ is a constant. Hence, $m+n\geq-
2\log_{2}x-N$. Since obviously $r_{m}\leq1$ we have that
$$
r_{m}(t,x)\leq 1\wedge
(N\gamma^{-m+M^{-}_{-2\log_{2}x-N}}).
$$
By the assumption if $x$ is small enough
$$
M^{-}_{-2\log_{2}x-N}>\alpha(-2\log_{2}x).
$$
Therefore, for $x\in(0,2^{-m/2}]$
\begin{equation}
                                             \label{2.25.5}
r_{m}(t,x)\leq 1\wedge
(N\gamma^{-m -2\alpha\log_{2}x  })=1\wedge(N\gamma^{-m}x^{\chi}).
\end{equation}

Next,
$$
\int_{0}^{2^{- m/2}}\frac{1}{x^{\nu p}}r^{p}_{m}(t,x)
\,dx
\leq \int_{0}^{\infty}
 \frac{1}{x^{\nu p}}(1\wedge(N\gamma^{-m}x^{\chi}))^{p}\,dx
$$
$$
=\gamma^{m(1-\nu p)/\chi}
\int_{0}^{\infty}
 \frac{1}{x^{\nu p}}(1\wedge(N x^{\chi}))^{p}\,dx,
$$
where the last integral is finite owing to \eqref{2.7.2}.
This proves the lemma.

Let $w_{t}$ be a Wiener process with respect to a
filtration $\{\cF_{t},t\geq0\}$  of complete
$\sigma$-fields and let $a_{t}$ and $\sigma_{t}$ be bounded
real-valued  processes
predictable with respect to $\{\cF_{t},t\geq0\}$  and such that
$a_{t}-\sigma^{2}_{t}\geq\delta\sigma^{2}_{t}$, 
where $\delta\in(0,\infty)$ is a constant, 
$a_{t}-\sigma^{2}_{t}>0$ for all $(\omega,t)$ and
for all $\omega$
$$
\int_{0}^{\infty}[a_{t}-\sigma^{2}_{t}]\,dt=\infty.
$$

Set $D_{x}=\partial/\partial x$.
For $m=0,1,2,...$ we will be dealing with
the SPDE
$$
dv(t,x)=(1/2)a_{t}D^{2}_{x }v(t,x)\,dt+\sigma_{t}
D_{x}v(t,x)\,dw_{t}
$$
in $B_{m}=(0,\infty)\times(0,2^{-m/2})$ with boundary conditions 
\begin{equation}
                                           \label{3.19.2}
v(t,0)=0,\quad v(t,2^{-m/2})=1,\quad t>0,
\end{equation}
\begin{equation}
                                           \label{3.19.3}
v(0,x)=0,\quad0<x<1.
\end{equation}

Recall that by Theorem 2.1 of \cite{Kr03}
there is a deterministic function $\alpha_{0}(c)$, $c>0$,
such that $\alpha_{0}(c)\to1$ as $c\to\infty$ and
with probability one for any $T\in(0,\infty)$
$$
\nliminf_{n\to\infty}\inf_{t\in[0,T]}
\frac{1}{n+1}M_{n}^{-}(w_{\cdot},c,t)=\alpha_{0}(c).
$$

\begin{theorem}
                                      \label{theorem 1.26.1}
For each $m=0,1,2,...$
there is a function $v_{m}(t,x)=v_{m}(\omega,t,x)$
defined on $\Omega\times\bar{B}_{m}$ such that

(i) $v_{m}(t,x)$ is $\cF_{t}$-measurable for 
each $(t,x)\in\bar{B}_{m}$,

(ii) $v_{m}(t,x)$ is bounded and continuous in $\bar{B}_{m}
\setminus\{(0,2^{-m/2})\}$ for each $\omega$,

(iii) derivatives of $v_{m}(t,x)$ of any order with respect to $x$
are continuous in $B_{m}\cup(\{0\}\times(0,2^{-m/2}))$
 for each $\omega$,

(iv) equations  \eqref{3.19.2}  and  \eqref{3.19.3}  hold for
each $\omega$,

(v) almost surely, for any $(t,x)\in B_{m}$
$$
v_{m}(t,x)=\int_{0}^{t}(1/2)a_{s}D^{2}_{x }
v_{m}(s,x)\,ds+\int_{0}^{t}\sigma_{s}
D_{x}v_{m}(s,x)\,dw_{s},
$$

(vi) for any $T\in(0,\infty)$,   $c,d>0$, $p>0$,
$\alpha>0$
such that $\alpha_{0}(c\sqrt{\delta})>\alpha$, 
and
$\nu$ satisfying
\begin{equation}
                                           \label{1.26.01}
1<\nu p<\chi p+1  ,
\end{equation}
where $\chi=-2\alpha\log_{2}\gamma(c,d,1)$,
we have that with probability one
\begin{equation}
                                           \label{1.26.1}
\pi_{T}:=\sup_{m\geq0}\sup_{t\in[0,T]}2^{- m(\nu p-1)/(2\alpha)}
\int_{0}^{2^{-  m/2}}
\frac{1}{x^{\nu p}}v^{p}_{m}(t,x)
\,dx<\infty.
\end{equation}
 
\end{theorem}

 Proof. In Lemma \ref{lemma 1.25.1}
take $\delta=1$ and
set $\tilde{v}_{m}(x_{\cdot},x_{t}+x,t)=
r_{m}(x_{\cdot},t,x)$, where $r_{m}$
is introduced in that lemma.
Set
$$
\psi_{t}=\int_{0}^{t}(a_{s}-\sigma^{2}_{s})\,ds,
\quad\xi_{t}=\int_{0}^{\phi_{t}}\sigma_{s}\,dw_{s},
\quad\tilde{\cF}_{t}=\cF_{\phi_{t}},
$$
$$
\bar{v}_{m}(t,x)=\bar{v}_{m}(\omega,t,x)
=\tilde{v}_{m}(\xi_{\cdot},t,x),
\quad 
v_{m}(t,x)=v_{m}(\omega,t,x) 
=\bar{v}_{m}(\psi_{t},x+\xi_{\psi_{t}}),
$$
where $\phi_{t}=\inf\{s\geq0:\psi_{s}\geq t\}$  
 is the inverse function to $\psi_{t}$.

It is proved in Theorem 4.1 of \cite{Kr03}
that $v_{0}$ possesses properties (i)-(v). The proof
that this is also true for any $m$
is no different.

Furthermore, it is well known that
$$
\sqrt{\delta}
\int_{0}^{t}\sigma_{s}\,dw_{s}=\tilde{w}_{\tilde{\psi}(t)},
$$
where $\tilde{w}_{t}$ is a Wiener process and
$$
\tilde{\psi}_{t}=\delta\int_{0}^{t}\sigma^{2}_{s}\,ds.
$$
Hence $\xi_{t}=\delta^{-1/2}\tilde{w}_{\tilde{\psi}(\phi_{t})}$
with 
$$
(\tilde{\psi}(\phi_{t}))'=\delta
\sigma^{2}_{s}/(a_{s}-\sigma^{2}_{s})
|_{s=\phi_{t}}\leq 1.
$$
It follows that  for $n=0,1,2,...$
we have
$$
M^{-}_{n}(\xi_{\cdot},c,t)\geq M^{-}_{n}
(\tilde{w}_{\cdot},c\sqrt{\delta},\tilde{\psi}(\phi_{t})),
$$
$$
\inf_{t\leq T}M^{-}_{n}(\xi_{\cdot},c,t)\geq 
\inf_{t\leq T} M^{-}_{n}
(\tilde{w}_{\cdot},c\sqrt{\delta},t) ,
$$
and with probability one
$$
\nliminf_{n\to\infty}\inf_{t\leq T}
\frac{1}{n+1}M^{-}_{n}(\xi_{\cdot},c,t)\geq\alpha_{0}
(c\sqrt{\delta})>\alpha.
$$

Finally, for $M=\sup_{\omega,t}(a_{t}-\sigma^{2}_{t})$
we have
$$
\sup_{t\leq T}\int_{0}^{2^{- m/2}}\frac{1}{x^{\nu p}}
v^{p}_{m}(t,x)\,dx
\leq
\sup_{t\leq MT}\int_{0}^{2^{- m/2}}\frac{1}{x^{\nu p}}
\bar{v}^{p}_{m}(t,x+\xi_{t})\,dx
$$
$$
=\sup_{t\leq MT}\int_{0}^{2^{-  m/2}}\frac{1}{x^{\nu p}}
\tilde{v}^{p}_{m}(\xi_{\cdot},t,x+\xi_{t})\,dx
=\sup_{t\leq MT}\int_{0}^{2^{-  m/2}}\frac{1}{x^{\nu p}}
r^{p}_{m}(\xi_{\cdot},t,x )\,dx.
$$
After this it only remains to use Lemma \ref{lemma 1.25.2}.
The theorem is proved.

\begin{remark}
                                         \label{remark 2.23.01}
Obviously, for any $\varepsilon\in(0,2^{-(m+2)/2})$
we have 
$$
v_{m}(t,\cdot)\in W^{1}_{2}(\varepsilon,2^{-m/2}
-\varepsilon)
$$
 for any $t\in[0,\infty)$
and for any $T\in(0,\infty)$ we have
$$
\int_{0}^{T}\|v_{m}(t,\cdot)\|_{W^{1}_{2}(\varepsilon,2^{-m/2}
-\varepsilon)}\,dt<\infty.
$$
Furthermore, by using the deterministic and stochastic
versions of Fubini's theorem one easily proves that
for any $\phi\in C^{\infty}_{0}(0,2^{-m/2})$ with probability one
for all $t\in[0,\infty)$
$$
(\phi,v_{m}(t,\cdot)=(1/2)\int_{0}^{t}(\phi,
a_{s}D^{2}_{x}v_{m}(s,\cdot))\,ds+
\int_{0}^{t}(\phi,\sigma_{s}D_{x}v_{m}(s,\cdot))\,dw_{s}.
$$
\end{remark}

\mysection{Continuity of solutions of SPDEs}
                                        \label{section 2.22.2}

We take the processes $a_{t}$, $\sigma_{t}$
as before Theorem \ref{theorem 1.26.1} but impose
stronger assumptions on them.

Assume that
there exist  constants $\delta_0,\delta_{1}\in(0,1]$ such that,
for every $(\omega,t)$  
$$
\delta_0 \leq\delta_{1}  a _{t} 
 \leq a _{t}- \sigma^{2}_{t} 
\leq \delta_0^{-1}  ,
$$
We will be dealing with solutions $u_{t}(x)$   of
\begin{equation}
                                          \label{2.22.3}
du_{t}=((1/2)a_{t}D_{x}^{2}u_{t}+f_{t})\,dt+
(\sigma_{t}D_{x}u_{t}+g_{t})
\,dw_{t}
\end{equation}
on $\bR_{+}$ with zero initial condition.
To specify the assumptions on $f,g$ and the class
of solutions we borrow the Banach spaces
$\bH^{\gamma}_{p,\theta}(\tau)$ and
$\bL _{p,\theta}(\tau)$ 
from \cite{KL98}. 
We also denote by $M$ the operator of multiplying by $x$.
Recall that, for $p\geq2,0<\theta<p$,
 the norms in $\bH^{\gamma}_{p,\theta}(\tau)$,
$\gamma=1,2$, and $\bL _{p,\theta}(\tau)$ are given by
$$
\| v\|^{p}_{\bL_{p,\theta}(\tau)}=E\int_{0}^{\tau}
 \int_{0}^{\infty}x^{\theta-1}|v(t,x) |^{p}\,dxdt,
$$
$$
\| v\|_{\bH^{1}_{p,\theta}(\tau)}=
\| v\|_{\bL_{p,\theta}(\tau)}+
\| MD_{x}v\|_{\bL_{p,\theta}(\tau)},
$$
$$
\| v\|_{\bH^{2}_{p,\theta}(\tau)}
=\| v\|_{\bH^{1}_{p,\theta}(\tau)}+
\| M^{2}D^{2}_{x}v\|_{\bL_{p,\theta}(\tau)}.
$$

Given  $p\geq2$,
  $\theta\in[p-1,p)$,
any stopping time $\tau$, $f\in M^{-1}\bL_{p,\theta}(\tau)$, 
and $g\in\bH^{1}_{p,
\theta}(\tau)$
by Theorem 3.2 of \cite{KL98} equation \eqref{2.22.3}
with zero initial condition
has a unique solution $u\in M\bH^{2}_{p,\theta}(\tau)$ and
$$
\|M^{-1}u\|_{\bH^{2}_{p,\theta}(\tau)}
\leq N(\|M f\|_{\bL _{p,\theta}(\tau)}+
\|g\|_{\bH^{1}_{p,\theta}(\tau)}),
$$
where $N=N(p,\theta,\delta_{0},\delta_{1})$.

We will also use Theorem 4.7 of \cite{Kr01},
which implies that if $u$ is a solution of
\eqref{2.22.3} of class $M\bH^{2}_{p,\theta}(\tau)$
with zero initial condition and $f\in M^{-1}\bL_{p,\theta}(\tau)$, 
and $g\in\bH^{1}_{p,
\theta}(\tau)$ and if there are numbers
$T\in(0,\infty)$ and $\beta$ such that
$$
2/p<\beta\leq1,\quad \tau\leq T,
$$
then for almost any $\omega$ the function
$u_{t}(x)$ is continuous in $(t,x)$ (that is, has a continuous
modification) and
$$
E\sup_{t\leq\tau}\sup_{x>0}|x^{\beta-1+\theta/p}
u_{t}(x)|^{p}\leq NT^{\beta p/2}(\|M^{-1}u\|_{\bH^{2}_{p,\theta}(\tau)}
+\|M f\|_{\bL _{p,\theta}(\tau)}+
\|g\|_{\bH^{1}_{p,\theta}(\tau)}),
$$
where $N=N(p,\theta,\beta,\delta_{0})$.

 Everywhere below we take 
$$
p>2.
$$
\begin{theorem}
                                           \label{theorem 2.7.1}
Let   $T\in(0,\infty)$,
 $c>0$, $\alpha\in(0,1)$, $\theta>0$, $\mu$ be some constants
such that $\alpha_{0}(c\sqrt{\delta_{1}})>\alpha$,
$$
\theta_{0}<\theta<p,\quad
\mu< p(1+2\log_{2}\gamma(c ))-2
=\theta_{0}-2+2p(1-\alpha)\log_{2}\gamma(c),
$$
where
$$
\gamma(c)=\gamma(c,1,1),\quad
\theta_{0}=p(1+2\alpha\log_{2}\gamma(c))\quad(>0).
$$
Let $f\in M^{-1}\bL_{p,p-1}(T)$, 
 $g\in\bH^{1}_{p,
p-1}(T)$, and let $u\in M\bH^{2}_{p,
p-1}(T)$ be a solution of \eqref{2.22.3}
with zero initial condition.
 
Finally, assume that $f_{t}(x)=g_{t}(x)=0$ for $x\geq1$
and $f\in M^{-1}\bL_{p,\mu}(T)$, 
 $g\in\bH^{1}_{p,\mu}(T)$. Then
  there exist    stopping times $\tau_{n}\uparrow T$,
defined independently of $f $ and $g $ such that,
for each $n$,
$u\in M\bH^{2}_{p,\theta}(\tau_{n})$ and
\begin{equation}
                                          \label{2.7.02}
\|M^{-1}u\|^{p}_{\bH^{2}_{p,\theta}(\tau_{n})}
\leq n(\|Mf\|^{p}_{\bL_{p,\mu}(\tau_{n})}+
\|g\|^{p}_{\bH^{1}_{p,\mu}(\tau_{n})})
\end{equation}

\end{theorem}

Here is the result about the continuity of $u_{t}(x)$
we were talking about in the introduction.
\begin{remark}
                                           \label{remark 2.23.1}
By Theorem 4.7 of \cite{Kr01} and Theorem \ref{theorem 2.7.1}
if we have a number $\beta\in(2/p,1]$, then
there exists a sequence of stopping times $\tau_{n}
\uparrow T$ such that
$$
E\sup_{t\leq \tau_{n}}\sup_{x>0}|x^{-\varepsilon}u_{t}(x)|^{p}
<\infty,
$$
where $\varepsilon= 1-\beta-\theta/p$. Due to the freedom
of choosing $\alpha$,  
$\beta$,  and $\theta$, the number $\varepsilon$
can be made as close from the right as we wish to
$$
1-\lim_{\alpha\to\alpha_{0}(c\sqrt{\delta_{1}})}
(2+\theta_{0})/p =-2/p-2\alpha_{0}(c\sqrt{\delta_{1}})
\log_{2}\gamma(c).
$$
If we allow arbitrary $p$, then the rate
of convergence of $u_{t}(x)$ to zero as $x\downarrow0$
is almost 
$$
x^{ \varepsilon_{0}},\quad
\varepsilon_{0}=-2\alpha
_{0}(c\sqrt{\delta_{1}})\log_{2}\gamma(c) >0,
$$
which is the same as we obtained for 
$v_{m}(t,x)$ (see \eqref{2.25.5}). Hence, the presence
of $f$ and $g$ does not spoil the situation too much.

It is also worth noting that $f$ and $g$ still may blow up
near zero even if $p$ is large. When $p$ is large
we can take $(\mu-1)/p$ as close to $1+2\log_{2}\gamma(c)$
as we wish and then the integral
$$
\int_{0}^{1}x^{\mu-1}|xf_{t}(x)|^{p}\,dx
$$
converges if $|f_{t}(x)|$
blows up near $x=0$ slightly slower than
$ x^{-2(1+\log_{2}\gamma(c))}$. Here $\log_{2}\gamma(c)\to0$
as $c\to\infty$ and one can allow $|f_{t}(x)|$ to blow up
almost as $x^{-2}$.

However, when $f$ and $g$ become more irregular near 0, the
rate with which the solution goes to zero at 0
deteriorates. In connection with this it is interesting
to investigate what happens with $\varepsilon_{0}$
as $\delta_{1}\downarrow0$. Take an $m$ so large that
$\alpha_{0}(m)>1/2$ and set $c=m\delta_{1}^{-1/2}-1/2$. Then
for $\delta_{1}$ small we have $
\alpha_{0}(c\sqrt{\delta_{1}})>1/2$ and
$$
\varepsilon_{0}\geq-\log_{2}[1-
P(\min_{s\leq1/2}w_{s}\leq-c-1/\sqrt{2},\max_{s\leq1/2}w_{s}\leq
1-1/\sqrt{2})),
$$
$$
\varepsilon_{0}\ln2\geq-\ln[1-
P(\min_{s\leq1/2}w_{s}\leq-c-1/\sqrt{2},\max_{s\leq1/2}w_{s}\leq
1-1/\sqrt{2}))
$$   
$$
\sim P(\min_{s\leq1/2}w_{t}\leq-c-1/\sqrt{2},
\max_{s\leq1/2}w_{s}\leq
1-1/\sqrt{2})
$$
$$
=P(\min_{s\leq1/2}w_{t}\leq-c-1/\sqrt{2})-
P(\min_{s\leq1/2}w_{t}\leq-c-1/\sqrt{2},
\max_{s\leq1/2}w_{s}\geq
1-1/\sqrt{2})
$$
and 
$$
P(\min_{s\leq1/2}w_{t}\leq-c-1/\sqrt{2},
\max_{s\leq1/2}w_{s}\geq
1-1/\sqrt{2})\leq2P(\min_{s\leq1/2}w_{t}\leq-c-1 ),
$$
so that
$$
 P(\min_{s\leq1/2}w_{t}\leq-c-1/\sqrt{2},
\max_{s\leq1/2}w_{s}\leq
1-1/\sqrt{2})
$$
$$
\geq 
P(\min_{s\leq1/2}w_{t}\leq-c-1/\sqrt{2})
-2P(\min_{s\leq1/2}w_{t}\leq-c-1 ).
$$
Next, as $a\to\infty$
$$
P(\min_{s\leq1/2}w_{s}\leq-a)=P(|w_{1/2}|\geq a)
=\frac{2}{\sqrt{\pi}}\int_{a}^{\infty}e^{-x^{2}}\,dx
\sim\frac{1}{\sqrt{\pi}}a^{-1}e^{-a^{2}}
$$
and
$$
\lim_{\delta_{1}\downarrow0}
[P(\min_{s\leq1/2}w_{t}\leq-c-1/\sqrt{2})
$$
$$
-2P(\min_{s\leq1/2}w_{t}\leq-c-1 )](c+1/\sqrt{2})
e^{(c+1/\sqrt{2})^{2}}=\frac{1}{\sqrt{\pi}}.
$$

Hence
$$
\nliminf_{\delta_{1}\downarrow0} 
[m\delta_{1}^{-1/2}e^{m^{2}/\delta_{1}}
\varepsilon_{0}]\geq\frac{1}{\sqrt{\pi}\ln2}.
$$

This result may seem unsatisfactory since the
guaranteed value of $\varepsilon_{0}$
is extremely small when $\delta_{1}$ is small.
However, recall that by Remark 4.2 of \cite{Kr03}
the best possible rate with which the solutions go to zero
for small $\delta_{1}$
is less than
$$
(1+\kappa)(2\pi\delta_{1})^{-1/2}e^{-1/(2\delta_{1})},
$$
where $\kappa>0$ is any number.

\end{remark}

To prove Theorem \ref{theorem 2.7.1}, first we prove the following.
\begin{lemma}
                                           \label{lemma 2.22.1}
Assume that, for an $m=0,1,2,...$
we have
 $f_{t}(x)=g_{t}(x)=0$ if $x\leq2^{-m/2}$.
  Then almost surely for all $t\leq T$
and  $x\in(0,2^{-m/2})$
\begin{equation}
                                               \label{2.22.1}
|u_{t}( x)|\leq v_{m}(t,x)\sup_{s\leq t}|u_{s}( 2^{-m/2})|.
\end{equation}
 
\end{lemma}

Proof. 
By Theorem 4.7 of \cite{Kr01} the function $u_{t}(x)$
is continuous in $[0,T]\times(0,2^{-m/2})$ (a.s.) and therefore
to prove \eqref{2.22.1} it suffices to prove
that for each $\varepsilon\in(0,2^{-(m+2)/2})$
 almost surely for all $t\leq T$
and  $x\in D:=(\varepsilon,2^{-m/2}-\varepsilon)$
\begin{equation}
                                               \label{2.22.2}
|u_{t}^{\varepsilon}(x)|\leq 
v_{m}(t,x) 
\sup_{s\leq t}|u^{\varepsilon}_{s}( 2^{-m/2}-\varepsilon)| 
=:v_{m}(t,x)\rho^{\varepsilon}_{t},
\end{equation}
where $u^{\varepsilon}_{t}(x)=u_{t}(x-\varepsilon)$.
The function $u^{\varepsilon}_{t}$ satisfies 
\eqref{2.22.3} with $f=g=0$ in $(0,T)\times(\varepsilon,
2^{-m/2}+\varepsilon)$ and in 
 $(0,T)\times D$. 
Furthermore, (a.s.) for almost any $t\in(0,T)$
we have $D_{x}u_{t}\in L_{p}(D)$ implying that
the limit of $u_{t}^{\varepsilon}(x)$ as $x\downarrow\varepsilon$
exists. Since (a.s.) for almost all $t\in(0,T)$
also $(x-\varepsilon)^{-1}u_{t}^{\varepsilon}\in L_{p}(D)$,
the limit is zero. As $x\uparrow2^{-m/2}-\varepsilon$
the situation is simpler and we see that (a.s.)
for almost all $t\in(0,T)$ we have
$$
\lim_{D\ni x\to\partial D}
(u^{\varepsilon}_{t}(x)-v_{m}(t,x)\rho^{\varepsilon}_{t})^{+}=0.
$$ 

Furthermore, (a.s.) for almost all $t\in(0,T)$ it holds that
$u^{\varepsilon}_{t}\in W^{1}_{2}(D)$ and
$$
\int_{0}^{T}\|u^{\varepsilon}_{t}\|^{2}_{W^{1}_{2}(D)}\,dt<\infty.
$$
Combining this with Remark \ref{remark 2.23.01}
we see that (a.s.) for almost all $t\in(0,T)$ we have
$(u^{\varepsilon}_{t}-v_{m}(t,\cdot)
\rho^{\varepsilon}_{t})^{+}\in W^{1}_{2}(D)$,
$(u^{\varepsilon}_{t}-v_{m}(t,\cdot)
\rho^{\varepsilon}_{t})^{+}\in \WO{1}{2}(D)$ and
$$
\int_{0}^{T}\|
(u^{\varepsilon}_{t}-v_{m}(t,\cdot)\rho^{\varepsilon}_{t})^{+}
\|^{2}_{W^{1}_{2}(D)}\,dt<\infty.
$$
By Theorem \ref{theorem 3.18.1}  we conclude that
 almost surely for all $t\leq T$
and  $x\in D$
$$
u_{t}^{\varepsilon}(x)\leq v_{m}(t,x)\rho^{\varepsilon}_{t}.
$$
By combining this with similar
inequality for $-u_{t}^{\varepsilon}$ we obtain
\eqref{2.22.2}. The lemma is proved.

{\bf Proof of Theorem \ref{theorem 2.7.1}}. 
Clearly, we only need prove Theorem \ref{theorem 2.7.1}
for $f$ and $g$ such that $f_{t}(x)=g_{t}(x)=0$
for all $\omega,t$ if $x$ is small. Then
\begin{equation}
                                             \label{2.25.1}
f\in M^{-1}\bL_{p,\vartheta}(T),\quad
g\in\bH^{1}_{p,\vartheta}(T)
\end{equation}
for any $\vartheta$. 

According to Lemma 3.6 of \cite{KL98},
for each stopping time $\tau_{n}\leq T$,  
we have $u\in M\bH^{2}_{p,\theta}(\tau_{n})$ if
$u\in M\bL_{p,\theta}(\tau_{n})$
and under this condition 
the left-hand side of \eqref{2.7.02} is dominated
by a constant $N=N(\theta,p,\delta_{0},\delta_{1})$
times
\begin{equation}
                                             \label{2.25.2}
\|M^{-1}u\|^{p}_{\bL_{p,\theta}(\tau_{n})}
+E\int_{0}^{\tau_{n}}\int_{0}^{\infty}x^{\theta-1} 
|F_{t}(x) |^{p}  \,dxdt,
\end{equation}
where
$$
F_{t}(x):=
|xf_{t}(x) | +|g_{t}(x)| 
+|xD_{x}g_{t}(x)|.
$$

Observe that
obviously ($\alpha,
\gamma(c)\leq1$)
\begin{equation}
                                               \label{2.8.1}
\theta>\theta_{0}>\mu
\end{equation}
and since $f_{t}(x)=g_{t}(x)=0$ for $x\geq1$,
the integral involving $F_{t}$   will increase
if   we replace $\theta$ with $\mu$. It follows that
to prove the theorem,
it suffices
to estimate only the lowest norm of $u$, that is
to prove the existence of $\tau_{n}\uparrow T$ such that
\begin{equation}
                                          \label{2.7.002}
\|M^{-1}u\|^{p}_{\bL_{p,\theta}(\tau_{n})}\leq n
[\|Mf\|^{p}_{\bL_{p,\mu}(\tau_{n})}
+\|g\|^{p}_{\bH^{1}_{p,\mu}(\tau_{n})}].
\end{equation}

Next, take a $\vartheta\in[p-1,p)$ such that
$\vartheta>\theta$. For any stopping time $\tau\leq T$, by 
Lemma 4.3 of \cite{KL98} we have $u\in M\bH^{2}_{p,\vartheta}(\tau)$
and by Theorem 3.2 of  \cite{KL98}
\begin{equation}
                                               \label{2.8.3}
E\int_{0}^{\tau }\int_{0}^{\infty}x^{\vartheta-1}
|u_{t}(x)/x|^{p}
\,dxdt\leq N 
[\|Mf\|^{p}_{\bL_{p,\vartheta}(\tau )}
+\|g\|^{p}_{\bH^{1}_{p,\vartheta}(\tau )}],
\end{equation}
where $N=N(p,\vartheta,\delta_{0},\delta_{1})$. As before
on the right we can replace $\vartheta$ with $\mu$.
On the left one can replace $\vartheta$ with
$\theta$ if one restricts the domain of integration
with respect to $x$ to $x\geq1$.
Therefore \eqref{2.7.002} will be proved if  we
prove the existence of appropriate stopping times $\tau_{n}$
such that
\begin{equation}
                                          \label{2.7.3}
E\int_{0}^{\tau_{n}}\int_{0}^{1}x^{\theta-1}
|u_{t}(x)/x|^{p}
\,dxdt\leq n
[\|Mf\|^{p}_{\bL_{p,\mu}(\tau_{n})}
+\|g\|^{p}_{\bH^{1}_{p,\mu}(\tau_{n})}].
\end{equation}

Take a nonnegative
 $\eta\in C^{\infty}_{0}(\bR_{+})$ with support in
$(1,4)$ such that the $(1/2)$-periodic function
on $\bR$
$$
\sum_{k=-\infty}^{\infty}\eta(2^{x+k/2} )
$$
is identically equal to one. Introduce,
$$
\eta_{m}(x)=\eta(2^{m/2}x),\quad
(f_{mt},g_{mt}) =(f_{t} ,g_{t} )\eta_{m} .
$$

 Also   introduce $u_{mt}$
as solutions of class $M\bH^{2}_{p,p-1}(T)$ of \eqref{2.22.3} 
with zero initial condition and
  $f_{mt}$
and $g_{mt}$ in place of $f_{t}$ and $g_{t} $,
respectively.
Since only finitely many $f_{mt}$ and
$g_{mt}$ are not zero, we have
$$
u_{t}(x)=\sum_{m=1}^{\infty}u_{mt}(x)=I_{1}(t,x)+I_{2}(t,x),
$$
where
$$
I_{1}(t,x):=\sum_{m=1}^{\infty}u_{mt}(x)
I_{x\leq 2^{- m/2}},\quad I_{2}(t,x):=\sum_{m=1}^{\infty}u_{mt}(x)
I_{x>2^{- m/2}}.
$$

{\em Estimating $I_{2}$\/}. 
Take a $\vartheta$ as above, set 
$\varepsilon=(\vartheta-\theta)/(2p)$
and use H\"older's inequality to obtain
$$
|I_{2}(t,x)|^{p}\leq
\sum_{m=1}^{\infty}2^{\varepsilon pm}u_{mt}^{p}(x)J^{p/q}(x),
$$
where 
$$
J(x):=\sum_{m=1}^{\infty}
2^{-\varepsilon qm}I_{x>2^{- m/2}}\leq
Nx^{2\varepsilon q },\quad
J^{p/q}(x)\leq
Nx^{\vartheta-\theta} .
$$
Then use \eqref{2.8.3} again to get
$$
 E\int_{0}^{\tau}\int_{0}^{1} x^{\theta-1} |I_{2}(t,x)/x|^{p}
\,dxdt 
$$
$$
\leq N\sum_{m=1}^{\infty}
 E\int_{0}^{\tau}\int_{0}^{\infty}
2^{m(\vartheta-\theta)/2 }x^{\vartheta -1}
 |u_{mt}(x)/x|^{p}
\,dxdt 
$$
$$
\leq N\sum_{m=1}^{\infty}
 E\int_{0}^{\tau}\int_{0}^{\infty}
2^{m(\vartheta-\theta)/2 }x^{\vartheta -1}
 |F_{mt}(x)|^{p}
\,dxdt ,
$$
where
$$
F_{mt}(x)=|xf_{mt}(x)|+|g_{mt}(x)|+|xD_{x}g_{mt}(x)|.
$$
Here we notice few facts,
which will be also used in the future, that on the supports of
$f_{mt}(x)$ and $g_{mt}(x)$ we have $x\sim2^{-m/2}$,
$2^{m(\vartheta-\theta)/2 }F_{mt}(x)\sim
x^{\theta-\vartheta}F_{mt}(x)$ and
$$
F_{mt}(x)\leq F_{ t}(x) \bar{\eta}_{m}(x)
$$ 
where $\bar{\eta}_{m}(x)=
 \eta_{m}(x)+x2^{m/2}|\eta'(2^{m/2}x)|$. Notice that
the $(1/2)$-periodic function
$$
\sum_{m=-\infty}^{\infty}\bar{\eta}^{p}_{m}(2^{y})
$$
is bounded on $\bR$. Then we see that
$$
 E\int_{0}^{\tau}\int_{0}^{1} x^{\theta-1} |I_{2}(t,x)/x|^{p}
\,dxdt 
$$
$$
\leq N
 E\int_{0}^{\tau}\int_{0}^{\infty}
  x^{\theta -1}
|F_{ t}(x)|^{p} \sum_{m=1}^{\infty}\bar{\eta}^{p}_{m}(x)
\,dxdt
$$
$$
\leq N[\|Mf\|^{p}_{\bL_{p,\theta}(\tau )}
+\|g\|^{p}_{\bH^{1}_{p,\theta}(\tau )}]
$$
for any $\tau\leq T$ with a constant $N$ under control.
As above we can reduce $\theta$ in the last expression
to $\mu$.

{\em Estimating $I_{1}$\/}. Here we will see
how $\tau_{n}$ appear and how we get a
substantial drop
from $\theta$ to $\mu$. We have seen above that
the smaller $\mu$ is the weaker the statement
of the theorem becomes. Therefore, we may
concentrate on $\mu$ so close
to $p(1+2\log_{2}\gamma(c))-2$ from below that
$$
2<\beta p:=p(1+2\log_{2}\gamma(c))-\mu\leq p.
$$
Then  
$$
2/p<\beta\leq1.
$$

Observe that
$$
|I_{1}(t,x)|^{p}\leq\bigg(\sum_{m=1}^{\infty}
m^{-q}\bigg)^{p/q}\sum_{m=1}^{\infty}
m^{p}|u_{mt}(x)|^{p}I_{x\leq2^{- m/2}}
$$
$$
\leq N|\log_{2}x|^{p}\sum_{m=1}^{\infty}
 |u_{mt}(x)|^{p}I_{x\leq2^{- m/2}}.
$$
It follows that for any $\theta'<\theta$
$$
E\int_{0}^{\tau}\int_{0}^{1}x^{\theta-1}
|I_{1}(t,x)/x|^{p}\,dxdt
$$
$$
\leq N\sum_{m=1}^{\infty}E\int_{0}^{\tau}\int_{0}^{2^{-m/2}}
|\log_{2}x|^{p}x^{\theta-1}
|u_{mt}(x)/x|^{p}\,dxdt\leq NJ(\tau),
$$
where
$$
J(\tau):=\sum_{m=1}^{\infty}E\int_{0}^{\tau}\int_{0}^{2^{-m/2}}
 x^{\theta'-1}
|u_{mt}(x)/x|^{p}\,dxdt.
$$

By Theorem 4.7 of \cite{Kr01} and 
Theorem 3.3 of \cite{KL99}, for any $\tau\leq T$
\begin{equation}
                                                \label{2.24.3}
E\sup_{t\leq \tau}\sup_{x>0}|x^{\varepsilon}u_{mt}(x)|^{p}
\leq NT^{\beta p/2}E\int_{0}^{\tau}\int_{0}^{\infty}
x^{\theta_{1}-1}|F_{mt}(x)|^{p}\,dxdt,
\end{equation}
where $N=N(p,\delta_{0},\delta_{1},\beta)$ and
$$
\theta_{1}:=p-1,\quad \varepsilon:=\beta-1+\theta_{1}/p=\beta-1/p>0.
$$
 
Therefore,
$$
E\sup_{t\leq \tau}|u_{mt}(2^{- m/2})|^{p}
\leq2^{  m(\beta p-p+\theta_{1})/2}
E\sup_{t\leq \tau}\sup_{x>0}|x^{\varepsilon}u_{mt}(x)|^{p}
$$
$$
\leq N E\int_{0}^{\tau}\int_{0}^{\infty}
2^{  m(\beta p-p+\theta_{1})/2}
x^{\theta_{1}-1}|F_{mt}(x)|^{p}\,dxdt
$$
$$
\leq N 2^{m(\beta p-p+1)/2}E\int_{0}^{\tau}
\int_{0}^{\infty} 
|F_{mt}(x)|^{p}\,dxdt. 
$$

Next, observe that, by Lemma \ref{lemma 2.22.1}
for $x\in[0,2^{-m/2}]$ and $t\leq T$,
$$
|u_{mt}(x)|\leq v_{m}(t,x)\sup_{s\leq t}|u_{ms}(2^{-m/2})|.
$$
Hence
$$
J(\tau)\leq
\sum_{m=1}^{\infty}E\sup_{t\leq \tau}|u_{mt}( 2^{-m/2})|
\int_{0}^{\tau}\int_{0}^{2^{-m/2}}
 x^{\theta'-1}
|v_{m}(t,x)/x|^{p}\,dxdt
$$
$$
\leq
\sum_{m=1}^{\infty}2^{m(\nu p-1)/(2\alpha)}
E\sup_{t\leq \tau}|u_{mt}( 2^{-m/2})|
\int_{0}^{\tau}\pi_{t}\,dt,
$$
where $\nu $ is defined according to
$$
\nu p=p-\theta'+1
$$
and $\pi_{t}$ is introduced in Theorem  \ref{theorem 1.26.1}.
So far $\theta'$ was only restricted to $\theta'<\theta$,
so that $\nu p>1$. Due to the assumption that
$\theta>\theta_{0}$ one can satisfy $\theta_{0}
<\theta'<\theta$
in which case \eqref{1.26.01} holds. Then in light of Theorem
\ref{theorem 1.26.1}
  one can find   stopping times
$\tau_{n}\uparrow T$ such that
$$
\int_{0}^{\tau_{n}}\pi_{t}\,dt\leq n.
$$
Then
$$
J(\tau_{n})\leq nN \sum_{m=1}^{\infty}
2^{m(\nu p-1)/(2\alpha)}2^{m(\beta p-p+1)/2}
E
\int_{0}^{\tau_{n}}\int_{0}^{\infty} |F_{mt}(x)|
^{p}\,dxdt.
$$
As is easy to see the inequalities $\theta'>\theta_{0}$
and
$$
\beta p-p +(\nu p-1)/\alpha<-\mu 
$$
are equivalent.
Hence,
$$
E\int_{0}^{\tau_{n}}\int_{0}^{1}x^{\theta-1}
|I_{1}(t,x)/x|^{p}\,dxdt\leq NJ(\tau_{n})
$$
$$
\leq nN\sum_{m=1}^{\infty}2^{m(1-\mu)/2}
E
\int_{0}^{\tau_{n}}\int_{0}^{\infty} |F_{mt}(x)|
^{p}\,dxdt
$$
$$
\leq nN\sum_{m=1}^{\infty}E
\int_{0}^{\tau_{n}}\int_{0}^{\infty} x^{\mu-1}|F_{mt}(x)|
^{p}\,dxdt
$$
$$
\leq nN E
\int_{0}^{\tau_{n}}\int_{0}^{\infty} x^{\mu-1}
  |F_{ t}(x)|^{p} \,dxdt.
$$

By combining this estimate with the estimate of $I_{2}$,
noticing that the above constants $N$ are independent of
$f$ and $g$ and, if necessary,
 renumbering the sequence $\tau_{n}$ we come to
\eqref{2.7.3}. This proves the theorem.

\end{document}